# LATE POINTS FOR RANDOM WALKS
# IN TWO DIMENSIONS

By Amir Dembo,[1] Yuval Peres,[2] Jay Rosen[3]
and Ofer Zeitouni[4]

*Stanford University, University of California–Berkeley, City University of New York–College of Staten Island and Technion and University of Minnesota*

Let $\mathcal{T}_n(x)$ denote the time of first visit of a point $x$ on the lattice torus $\mathbb{Z}_n^2 = \mathbb{Z}^2/n\mathbb{Z}^2$ by the simple random walk. The size of the set of $\alpha$, $n$-late points $\mathcal{L}_n(\alpha) = \{x \in \mathbb{Z}_n^2 : \mathcal{T}_n(x) \geq \alpha \frac{4}{\pi}(n \log n)^2\}$ is approximately $n^{2(1-\alpha)}$, for $\alpha \in (0,1)$ [$\mathcal{L}_n(\alpha)$ is empty if $\alpha > 1$ and $n$ is large enough]. These sets have interesting clustering and fractal properties: we show that for $\beta \in (0,1)$, a disc of radius $n^\beta$ centered at nonrandom $x$ typically contains about $n^{2\beta(1-\alpha/\beta^2)}$ points from $\mathcal{L}_n(\alpha)$ (and is empty if $\beta < \sqrt{\alpha}$), whereas choosing the center $x$ of the disc uniformly in $\mathcal{L}_n(\alpha)$ boosts the typical number of $\alpha$, $n$-late points in it to $n^{2\beta(1-\alpha)}$. We also estimate the *typical* number of pairs of $\alpha$, $n$-late points within distance $n^\beta$ of each other; this typical number can be significantly smaller than the *expected* number of such pairs, calculated by Brummelhuis and Hilhorst [*Phys. A* **176** (1991) 387–408]. On the other hand, our results show that the number of ordered pairs of late points within distance $n^\beta$ of each other is larger than what one might predict by multiplying the total number of late points, by the number of late points in a disc of radius $n^\beta$ centered at a typical late point.

**1. Introduction.** Consider a simple random walk (SRW) on an $n \times n$ square with periodic boundary conditions (also called a lattice torus), run

Received April 2003; revised July 2004.

[1]Supported in part by NSF Grant DMS-00-72331 and by a US–Israel BSF grant.

[2]Supported in part by NSF Grant DMS-98-03597.

[3]Supported in part by grants from the NSF, from PSC-CUNY and from a US–Israel BSF grant.

[4]Supported in part by a US–Israel BSF grant.

until the "cover time," when it has visited every point of the square. Our focus will be on the set of uncovered points shortly before coverage, which we call "late points." In an important paper, Brummelhuis and Hilhorst [1] pointed out that in two dimensions, this set has an interesting fractal structure. The main finding of the present paper is that the set of late points has an even more subtle fractal structure than that suggested in [1]. A significant reason for this is that a key random variable measuring the structure of late points, namely the number of pairs of late points within distance $n^\beta$ of each other, has a median and mean of different orders of magnitude.

As noted in [1] this fractal structure is not present in three or higher dimensions, where at the scale of power laws the set of uncovered points resembles a uniformly sampled random set of the same size.

We proceed to a more quantitative discussion. Consider the SRW on the lattice torus $\mathbb{Z}_n^2 = \mathbb{Z}^2/n\mathbb{Z}^2$ starting at the origin. If $x \in \mathbb{Z}_n^2$, we let $\mathcal{T}_n(x)$ denote the time it takes the walk to first visit $x$. Let $\mathcal{T}_n = \max_{x \in \mathbb{Z}_n^2} \mathcal{T}_n(x)$ denote the time it takes the walk to completely cover $\mathbb{Z}_n^2$. In [4], Theorem 1.1, we showed that

$$(1.1) \qquad \lim_{n\to\infty} \frac{\mathcal{T}_n}{(n\log n)^2} = \frac{4}{\pi} \qquad \text{in probability.}$$

(Contrast this with the typical hitting time of a *fixed* point $x \in \mathbb{Z}_n^2$, which is of order $n^2 \log n$.)

We say that $x \in \mathbb{Z}_n^2$ is $\alpha, n$-late for some $0 < \alpha < 1$ if

$$\mathcal{T}_n(x) \geq \alpha \frac{4}{\pi}(n\log n)^2,$$

and set $\mathcal{L}_n(\alpha)$ to be the set of $\alpha, n$-late points in $\mathbb{Z}_n^2$. An adaptation of the arguments in [4] reveals that $|\mathcal{L}_n(\alpha)| \approx n^{2-2\alpha}$ in the following sense.

PROPOSITION 1.1. *For any $0 < \alpha < 1$,*

$$(1.2) \qquad \lim_{n\to\infty} \frac{\log |\mathcal{L}_n(\alpha)|}{\log n} = 2(1-\alpha) \qquad \text{in probability.}$$

If $\mathcal{L}_n(\alpha)$ were spread out uniformly in $\mathbb{Z}_n^2$, one would expect that for any $x \in \mathbb{Z}_n^2$ and $\alpha < \beta < 1$ we would have $|\mathcal{L}_n(\alpha) \cap D(x, n^\beta)| \approx n^{2\beta - 2\alpha}$. The next two theorems make precise the idea that the set $\mathcal{L}_n(\alpha)$ does not look like an independent uniform drawing of $n^{2-2\alpha}$ points in $\mathbb{Z}_n^2$, in the sense that $|\mathcal{L}_n(\alpha) \cap D(x, n^\beta)| \approx n^{2\beta - 2\alpha/\beta}$ for a typical $x$, whereas it is $\approx n^{2\beta(1-\alpha)}$ for most $x \in \mathcal{L}_n(\alpha)$.

THEOREM 1.2. *For any $0 < \alpha < \beta^2 < 1$ and $\delta > 0$,*

$$(1.3) \quad \lim_{n \to \infty} \max_{x \in \mathbb{Z}_n^2} \mathbf{P}\left( \left| \frac{\log |\mathcal{L}_n(\alpha) \cap D(x, n^\beta)|}{\log n} - (2\beta - 2\alpha/\beta) \right| > \delta \right) = 0.$$

In particular, for any $0 < \alpha, \beta < 1$ and any nonrandom sequence $x_n \in \mathbb{Z}_n^2$

$$(1.4) \quad \lim_{n \to \infty} \frac{\log |\mathcal{L}_n(\alpha) \cap D(x_n, n^\beta)|}{\log n} = \max(2\beta - 2\alpha/\beta, 0) \qquad \text{in probability.}$$

As stated already, the fractal nature of $|\mathcal{L}_n(\alpha)|$ is described by the next theorem that shows the clustering of late points; in the neighborhood of a "typical" $\alpha, n$-late point there is an "unusually large" number of $\alpha, n$-late points.

THEOREM 1.3. *For any $0 < \alpha, \beta < 1$ and $\delta > 0$,*

$$(1.5) \quad \lim_{n \to \infty} \max_{x \in \mathbb{Z}_n^2 \setminus \{0\}} \mathbf{P}\left( \left| \frac{\log |\mathcal{L}_n(\alpha) \cap D(x, n^\beta)|}{\log n} \right.\right.$$
$$\left.\left. - 2\beta(1 - \alpha) \right| > \delta \,\Big|\, x \in \mathcal{L}_n(\alpha) \right) = 0.$$

*Further, choosing $Y_n$ uniformly in $\mathcal{L}_n(\alpha)$,*

$$(1.6) \quad \lim_{n \to \infty} \frac{\log |\mathcal{L}_n(\alpha) \cap D(Y_n, n^\beta)|}{\log n} = 2\beta(1 - \alpha) \qquad \text{in probability.}$$

The predictions of [1], which motivated our work, are related to another description of the clustering properties of $\mathcal{L}_n(\alpha)$, obtained by focusing on pairs of late points.

THEOREM 1.4. *Let $0 < \alpha, \beta < 1$. Then*

$$(1.7) \quad \lim_{n \to \infty} \frac{\log |\{(x, y) \in \mathcal{L}_n^2(\alpha) : d(x, y) \leq n^\beta\}|}{\log n} = \rho(\alpha, \beta) \qquad \text{in probability,}$$

*where*

$$(1.8) \quad \rho(\alpha, \beta) = \begin{cases} 2 + 2\beta - 4\alpha/(2 - \beta), & \text{if } \beta \leq 2(1 - \sqrt{\alpha}), \\ 8(1 - \sqrt{\alpha}) - 4(1 - \sqrt{\alpha})^2/\beta, & \text{if } \beta \geq 2(1 - \sqrt{\alpha}). \end{cases}$$

For the mean number of pairs of $\alpha, n$-late points within distance $n^\beta$ of each other, Brummelhuis and Hilhorst ([1], (3.36)) obtain different growth exponents

$$(1.9) \quad \widehat{\rho}(\alpha, \beta) = \begin{cases} 2 + 2\beta - 4\alpha/(2 - \beta), & \text{if } \beta \leq 2 - \sqrt{2\alpha}, \\ 6 - 4\sqrt{2\alpha}, & \text{if } \beta \geq 2 - \sqrt{2\alpha}. \end{cases}$$

As we explain below, the functions

$$(1.10) \qquad F_{h,\beta}(\gamma) = \frac{(1-\gamma\beta)^2}{1-\beta} + h\gamma^2\beta,$$

of $\gamma \geq 0$, with $h$ a nonnegative integer, play an important role in the study of late points. It can be easily checked that

$$(1.11) \qquad \rho(\alpha,\beta) = 2 + 2\beta - 2\alpha \inf_{\gamma \in \Gamma_{\alpha,\beta}} F_{2,\beta}(\gamma),$$

where

$$(1.12) \qquad \Gamma_{\alpha,\beta} = \{\gamma \geq 0 : 2 - 2\beta - 2\alpha F_{0,\beta}(\gamma) \geq 0\}$$

(see Section 9). It is also easy to verify that

$$(1.13) \qquad \widehat{\rho}(\alpha,\beta) = \sup_{\beta' \leq \beta} \sup_{\gamma \geq 0} \{2 + 2\beta' - 2\alpha F_{2,\beta'}(\gamma)\},$$

so the difference between $\widehat{\rho}(\alpha,\beta)$ and $\rho(\alpha,\beta)$ is that the supremum in (1.13) is not subject to the constraint that $\gamma \in \Gamma_{\alpha,\beta}$. As explained below, this constraint differentiates the median number of pairs of $\alpha, n$-late points within distance $n^\beta$ of each other, easily obtained from (1.7), from its mean (found already in [1]).

The key to our approach lies in the following heuristic picture relating the lateness property to certain excursion counts for the random walk: fix an appropriate sequence of increasing radii $r_k$, $k = 1, \ldots, k_n$, with $r_{k+1}/r_k \sim r_k/r_{k-1}$, $r_0 = 1$ and $r_{k_n} \ll n$, and count the number of excursions $N_x(k)$ between $D(x, r_{k-1})$ and $D(x, r_k)$. A point that has many fewer than the typical number of excursions between these levels, by time $4\alpha(n\log n)^2/\pi$, is also extremely likely to be $\alpha$, $n$-late (see Lemma 4.1). Further, a *typical* $x \in \mathcal{L}_n(\alpha)$ has an atypical profile of excursion counts, determined approximately by considering a one-dimensional simple random walk on the set $\{1, \ldots, k_n\}$, started at $k_n$, and conditioned not to hit 1. Thus, not only is the point $x$ not hit by the random walk, but in fact a neighborhood of it is visited less often than it would have been otherwise, and this creates a large cluster of $\alpha, n$-late points in a neighborhood of such $x$.

Large deviations estimates for this one-dimensional walk imply that certain $\alpha$, $n$-late points $x$ have a much smaller number of excursions $N_x(\bar{k}_n)$ between discs in an intermediate scale $\bar{k}_n$, forcing an accumulation of many $\alpha, n$-late points in $D(x, r_{\bar{k}_n})$. In more detail, for $r_{\bar{k}_n} \approx n^\beta$, the probability of $N_x(\bar{k}_n)$ being near the value typically associated with $\alpha\gamma^2, n$-late points is about $n^{-2\alpha F_{0,\beta}(\gamma)}$. Given such a value of $N_x(\bar{k}_n)$, the probability that $x$ is an $\alpha, n$-late point is about $n^{-2\alpha\gamma^2\beta}$. Consequently, the probability of $x$ being $\alpha, n$-late with $N_x(\bar{k}_n)$ near the value typically associated with an $\alpha\gamma^2, n$-late point is about $n^{-2\alpha F_{0,\beta}(\gamma)}n^{-2\alpha\gamma^2\beta} = n^{-2\alpha F_{1,\beta}(\gamma)}$, and if we require that also a specific $y$ of distance $\approx n^\beta$ from $x$ is $\alpha$, $n$-late, the probability is further reduced to about $n^{-2\alpha F_{1,\beta}(\gamma)} n^{-2\alpha\gamma^2\beta} = n^{-2\alpha F_{2,\beta}(\gamma)}$. The constraint $\gamma \in \Gamma_{\alpha,\beta}$ in (1.11), which is missing in (1.13), represents the range of values of $N_x(\bar{k}_n)$ possibly found when examining all $O(n^{2-2\beta})$ centers $x$ of discs of radius $n^\beta$ that cover the torus $\mathbb{Z}_n^2$. Indeed, due to this constraint, the median of number of pairs of $\alpha$, $n$-late points within distance $n^\beta$ of each other is about $n^{\rho(\alpha,\beta)}$, whereas the mean of this variable is of the different order of magnitude $n^{\widehat{\rho}(\alpha,\beta)}$.

The value of $\rho(\alpha,\beta)$ is obtained by taking $\gamma \in \Gamma_{\alpha,\beta}$ for which the probability of locating specific pairs of $\alpha$, $n$-late points is maximal. This value of $\gamma$ coincides with the unconstrained minimizer of $F_{2,\beta}(\cdot)$ if and only if $\beta \leq 2(1 - \sqrt{\alpha})$, thus explaining the jump of $d^2\rho/d\beta^2$ at $\beta = 2(1 - \sqrt{\alpha})$. It is never the same as the typical $\gamma = 1$ [i.e., the minimizer of $F_{1,\beta}(\cdot)$], which one finds in most discs of radius $n^\beta$ centered at $\alpha$, $n$-late points. Hence, $\gamma = 1$ controls the exponent of Theorem 1.3. In contrast, the exponent of Theorem 1.2 is controlled by $\gamma = 1/\beta$ [i.e., the minimizer of $F_{0,\beta}(\cdot)$], found in most of the $O(n^{2-2\beta})$ discs of radius $n^\beta$ that cover $\mathbb{Z}_n^2$.

*Organization.* After a short section which collects some facts about the SRW, our paper is divided into three parts. The first part is about "global" properties of the set of $\alpha$, $n$-late points. It consists of Sections 3–5, where, adapting the arguments of [4], Sections 2, 3, 6, 7, to the context of simple random walk, we prove Proposition 1.1 and lay the groundwork for all other results. The second part deals with clustering of late points. It starts with the large deviation probability bounds of the form $n^{-2\alpha F_{h,\beta}(\gamma)}$, given in Section 6, which are key to our upper bounds, and moves on to the proofs of Theorems 1.2 and 1.3. The third part of the paper deals with Theorem 1.4 about pairs of $\alpha$, $n$-late points. Applying the bounds of Section 6 we derive the upper bound in Section 9, where we also solve the variational problem (1.11), with the complementary lower bound derived in Section 10 by a refinement of the construction of Section 4. In the final Section 11 we describe possible extensions of our results. We note that the arguments in this paper are based on direct analysis of the random walk, rather than a strong approximation argument with Brownian motion.

**2. Random walk preliminaries.** Let $S_n$, $n \geq 0$, denote a simple random walk (SRW) in $\mathbb{Z}^2$ and let $X_n$, $n \geq 0$, denote SRW in $\mathbb{Z}_K^2$. In this section we collect some facts about $S_n$, $n \geq 0$, and $X_n$, $n \geq 0$. We adopt here and throughout the paper the:

CONVENTION. Throughout, a function $Z(x)$ is said to be $O(x)$ if $Z(x)/x$ is bounded, uniformly in all implicit geometry-related quantities (such as $K$).

That is, $Z(x) = O(x)$ if there exists a universal constant $C$ (not depending on $K$) such that $|Z(x)| \leq Cx$. Thus $x = O(x)$ but $Kx$ is *not* $O(x)$. A similar convention applies to the symbol $o(x)$.

Let $D(x, r) = \{y \in \mathbb{Z}^2 : |y - x| < r\}$ where $|z|$ denotes the Euclidean norm of $z$. For any set $A \subseteq \mathbb{Z}^2$ we let $\partial A = \{y \in \mathbb{Z}^2 : y \in A^c, \text{and } \inf_{x \in A} |y - x| = 1\}$ and $\overline{A} = A \cup \partial A$. For any set $B \subseteq \mathbb{Z}^2$ let $T_B = \inf\{i \geq 0 : S_i \in B\}$ and $T'_B = \inf\{i \geq 1 : S_i \in B\}$. For $x, y \in A$ define the truncated Green function

$$G_A(x, y) = \sum_{i=0}^{\infty} \mathbb{E}^x(S_i = y, i < T_{\partial A}).$$

We have the following result which is Proposition 1.6.7 of [7]. For any $x \in D(0, n)$

$$(2.1) \qquad \mathbf{P}^x(T_0 < T_{\partial D(0,n)}) = \frac{\log(n/|x|) + O(|x|^{-1} + (\log n)^{-1})}{\log n}$$

and

$$(2.2) \qquad G_{D(0,n)}(x, 0) = \frac{2}{\pi} \log\left(\frac{n}{|x|}\right) + O(|x|^{-1} + n^{-1}).$$

We next note formula (1.21) of [7]: Uniformly for $x \in D(0, n)$,

$$(2.3) \qquad n^2 - |x|^2 \leq \mathbb{E}^x(T_{\partial D(0,n)}) \leq (n+1)^2 - |x|^2.$$

We also have the result of Exercise 1.6.8 of [7]: Uniformly in $r < |x| < R$,

$$(2.4) \qquad \mathbf{P}^x(T_{\partial D(0,r)} < T_{\partial D(0,R)}) = \frac{\log(R/|x|) + O(r^{-1})}{\log(R/r)}.$$

Define the hitting distribution of the boundary of $A$ by

$$H_{\partial A}(x, y) = \mathbf{P}^x(S_{T_{\partial A}} = y).$$

We have the following Harnack inequality.

LEMMA 2.1. *Uniformly for $\delta < 1/2$, $x, x' \in D(0, \delta n)$ and $y \in \partial D(0, n)$,*

$$(2.5) \qquad H_{\partial D(0,n)}(x, y) = (1 + O(\delta) + O(n^{-1}))H_{\partial D(0,n)}(x', y).$$

*Furthermore, if $\delta' < \delta$ are such that*

$$\min_{x \in \partial D(0, \delta n)} \mathbf{P}^x(T_{\partial D(0,n)} < T_{\partial D(0, \delta' n)}) \geq 1/4,$$

*then uniformly in $x \in \partial D(0, \delta n)$ and $y \in \partial D(0, n)$,*

$$(2.6) \qquad \begin{aligned} &\mathbf{P}^x(S_{T_{\partial D(0,n)}} = y, T_{\partial D(0,n)} < T_{\partial D(0, \delta' n)}) \\ &\quad = (1 + O(\delta) + O(n^{-1}))\mathbf{P}^x(T_{\partial D(0,n)} < T_{\partial D(0, \delta' n)})H_{\partial D(0,n)}(x, y). \end{aligned}$$

PROOF. By Lemma 1.7.3 of [7], for any $y \in \partial D(0, n)$ and $\delta < 1/2$,

$$H_{\partial D(0,n)}(x, y) = \sum_{z \in \partial D(0,n/2)} \mathbf{P}^y(S_{T'_{\partial D(0,n/2) \cup \partial D(0,n)}} = z)G_{D(0,n)}(z, x).$$

But

$$G_{D(0,(1-\delta)n)}(z - x, 0) \le G_{D(0,n)}(z, x) \le G_{D(0,(1+\delta)n)}(z - x, 0)$$

and by (2.2), with $|z - x| = n(1/2 + O(\delta))$,

$$\begin{aligned} G_{D(0,(1\pm\delta)n)}(z - x, 0) &= \frac{2}{\pi}\log\left(\frac{(1 \pm \delta)n}{|z - x|}\right) + O(n^{-1}) \\ &= \frac{2}{\pi}\log\left(\frac{1 \pm \delta}{1/2 + O(\delta)}\right) + O(n^{-1}) \end{aligned}$$

and (2.5) now follows.

Turning to (2.6), we have

$$\begin{aligned} (2.7) \quad &\mathbf{P}^x(S_{T_{\partial D(0,n)}} = y, T_{\partial D(0,n)} < T_{\partial D(0,\delta'n)}) \\ &\quad = H_{\partial D(0,n)}(x, y) - \mathbf{P}^x(S_{T_{\partial D(0,n)}} = y, T_{\partial D(0,n)} > T_{\partial D(0,\delta'n)}). \end{aligned}$$

By the strong Markov property at $T_{\partial D(0,\delta'n)}$,

$$\begin{aligned} (2.8) \quad &\mathbf{P}^x(S_{T_{\partial D(0,n)}} = y, T_{\partial D(0,n)} > T_{\partial D(0,\delta'n)}) \\ &\quad = \mathbb{E}^x(H_{\partial D(0,n)}(S_{T_{\partial D(0,\delta'n)}}, y); T_{\partial D(0,n)} > T_{\partial D(0,\delta'n)}). \end{aligned}$$

Since $\partial D(0, \delta n)$ separates $\partial D(0, n)$ from $\partial D(0, \delta'n)$, by the strong Markov property and (2.5), uniformly in $w \in \partial D(0, \delta'n)$,

$$\begin{aligned} H_{\partial D(0,n)}(w, y) &= \mathbb{E}^w(H_{\partial D(0,n)}(S_{T_{\partial D(0,\delta n)}}, y)) \\ &= (1 + O(\delta) + O(n^{-1}))H_{\partial D(0,n)}(x, y). \end{aligned}$$

Substituting back into (2.8) we have

$$\begin{aligned} &\mathbf{P}^x(S_{T_{\partial D(0,n)}} = y, T_{\partial D(0,n)} > T_{\partial D(0,\delta'n)}) \\ &= (1 + O(\delta) + O(n^{-1}))\mathbf{P}^x(T_{\partial D(0,n)} > T_{\partial D(0,\delta'n)})H_{\partial D(0,n)}(x, y). \end{aligned}$$

Combining this with (2.7) and the assumptions of the lemma, used to control the error terms, we obtain (2.6) which completes the proof of the lemma. □

Combining the above with Lemma 1.7.4 of [7] we see that if $\mu_n$ denotes uniform measure on $\partial D(0, n)$, then for all $\delta < 1/2$ and some constants $0 < c = c(\delta) < C = C(\delta) < \infty$ we have that uniformly for $x \in D(0, \delta n)$,

$$(2.9) \quad c\mu_n(\cdot) \le H_{\partial D(0,n)}(x, \cdot) \le C\mu_n(\cdot).$$

Let $\widehat{H}_A(z,x) = \mathbf{P}^z(X_{\widehat{T}'_A} = x)$ be the hitting measure on $A \subseteq \mathbb{Z}_K^2$ by $X_n$ with $\widehat{T}_A$ and $\widehat{T}'_A$ the corresponding hitting times. When dealing with $X_n$, sets such as $D(x,r)$ and $\partial D(x,r)$ are defined with respect to the $L^2$-distance $d(\cdot, \cdot)$ in $\mathbb{Z}_K^2$.

LEMMA 2.2. *Uniformly in $K$, $z, z' \in \partial D(0, R)$ and $x \in \partial D(0, r)$ with $4r < R < K/2$,*

$$(2.10) \qquad \widehat{H}_{\partial D(0,r)}(z,x) = \left(1 + O\left(\frac{r}{R}\log\frac{R}{r}\right)\right)\widehat{H}_{\partial D(0,r)}(z',x).$$

*Furthermore, if $4r < R < R' < K/2$ are such that*

$$\min_{z \in \partial D(0,R)} \mathbf{P}^z(T'_{\partial D(0,r)} < T'_{\partial D(0,R')}) \geq 1/4,$$

*then uniformly in $z \in \partial D(0,R)$ and $x \in \partial D(0,r)$,*

$$
\begin{aligned}
(2.11) \quad & \mathbf{P}^z(X_{T'_{\partial D(0,r)}} = x; T'_{\partial D(0,r)} < T'_{\partial D(0,R')}) \\
& = \left(1 + O\left(\frac{r}{R}\log\frac{R}{r}\right)\right)\mathbf{P}^z(T'_{\partial D(0,r)} < T'_{\partial D(0,R')})\widehat{H}_{\partial D(0,r)}(z,x),
\end{aligned}
$$

*and if in addition $r^{-1} = O(\frac{r}{R})$, then uniformly in $z, z' \in \partial D(0, R)$ and $x \in \partial D(0,r)$,*

$$
\begin{aligned}
(2.12) \quad & \mathbf{P}^z(X_{T'_{\partial D(0,r)}} = x; T'_{\partial D(0,r)} < T'_{\partial D(0,R')}) \\
& = \left(1 + O\left(\frac{r}{R}\log\frac{R}{r}\right)\right)\mathbf{P}^{z'}(X_{T'_{\partial D(0,r)}} = x; T'_{\partial D(0,r)} < T'_{\partial D(0,R')}).
\end{aligned}
$$

PROOF. The bounds of (2.10) will follow immediately from the fact that uniformly in $z \in \partial D(0, R)$ and $x \in \partial D(0,r)$,

$$
\begin{aligned}
(2.13) \quad & \widehat{H}_{\partial D(0,r)}(z,x) = \left(1 + O\left(\frac{r}{R}\log\frac{R}{r}\right)\right) \\
& \qquad\qquad \times \frac{\mathbf{P}^x(T'_{\partial D(0,r)} > T'_{\partial D(0,R/2)})}{\sum_{x' \in \partial D(0,r)} \mathbf{P}^{x'}(T'_{\partial D(0,r)} > T'_{\partial D(0,R/2)})}.
\end{aligned}
$$

This is the equation above Theorem 2.1.3 of [7]. However, since that equation deals with the simple random walk in $\mathbb{Z}^2$ and $\widehat{H}_{\partial D(0,r)}(z,x)$ involves paths for which the difference between $\mathbb{Z}^2$ and $\mathbb{Z}_K^2$ might be significant, we next explain why the same proof works for $\mathbb{Z}_K^2$.

The proof of Lemma 2.1.1 of [7] shows that, with $A = \partial D(0, r)$, $B = \partial D(0, R/2)$ and $z \in \partial D(0, R)$,

$$\widehat{H}_A(z,x) = \frac{\sum_{v \in B} \widehat{G}_{\overline{D}(0,r)^c}(z,v)\widehat{H}_{A \cup B}(v,x)}{\sum_{v \in B} \widehat{G}_{\overline{D}(0,r)^c}(z,v)\mathbf{P}^v(\widehat{T}'_A < \widehat{T}'_B)},$$

with $\widehat{G}_{\overline{D}(0,r)^c}(z,v)$ the Green's function for $\overline{D}(0,r)^c$, the complement of $\overline{D}(0,r)$ in $\mathbb{Z}_K^2$. But this gives

$$(2.14) \qquad \min_{v \in B} \frac{\widehat{H}_{A \cup B}(v,x)}{\mathbf{P}^v(\widehat{T}'_A < \widehat{T}'_B)} \le \widehat{H}_A(z,x) \le \max_{v \in B} \frac{\widehat{H}_{A \cup B}(v,x)}{\mathbf{P}^v(\widehat{T}'_A < \widehat{T}'_B)}.$$

Note that $B = \partial D(0, R/2)$ separates $A = \partial D(0,r)$ from the complement of $\overline{D}(0, R/2)$ in $\mathbb{Z}_K^2$. Hence, the above max and min involve expressions that are determined by paths confined between $A = \partial D(0,r)$ and $B = \partial D(0, R/2)$, which are thus the same for the simple random walks in $\mathbb{Z}^2$ and in $\mathbb{Z}_K^2$. Consequently, (2.14) is precisely the top inequality on page 49 of [7], from which (2.13) follows. This completes the proof of (2.10). The bounds of (2.11) follow from (2.10) in the same way that (2.6) follows from (2.5). Finally, combining (2.10), (2.11) and (2.4) leads to (2.12). $\quad\square$

We next show that for $R' \gg R \gg r \gg 1$, the $\sigma$-algebra of excursions of the path from $\partial D(0,r)$ to $\partial D(0,R)$, prior to $T_{\partial D(0,R')}$, is almost independent of the initial point $z \in \partial D(0,R)$ and the final point $w \in \partial D(0,R')$.

LEMMA 2.3. *For $4r < R < R' < K/2$ and a random walk path starting at $z \in D(0,R)$, let $\mathcal{H}$ denote the $\sigma$-algebra generated by the excursions of the path from $\partial D(0,r)$ to $\partial D(0,R)$, prior to $T_{\partial D(0,R')}$. Suppose $r^{-1} = O(\frac{r}{R})$ and $\log(R'/R) \ge (1/4)\log(R/r)$. Then, uniformly in $K$, $z, z' \in \partial D(0,R)$, $w \in \partial D(0,R')$ and $B \in \mathcal{H}$,*

$$(2.15) \qquad \mathbf{P}^z(B | X_{T_{\partial D(0,R')}} = w) = \left(1 + O\left(\frac{R}{R'}\right)\right)\mathbf{P}^z(B)$$

*and*

$$(2.16) \qquad \mathbf{P}^z(B) = \left(1 + O\left(\frac{r}{R}\log\frac{R}{r}\right)\right)\mathbf{P}^{z'}(B).$$

PROOF. Fixing $z \in \partial D(0,R)$ it suffices to consider $B \in \mathcal{H}$ for which $\mathbf{P}^z(B) > 0$. Fix such $B$ and a point $w \in \partial D(0,R')$. Let $\tau_0 = 0$ and for $i = 0, 1, \ldots$ define

$$\tau_{2i+1} = \inf\{t \ge \tau_{2i} : S_t \in \partial D(0,r) \cup \partial D(0,R')\},$$

$$\tau_{2i+2} = \inf\{t \ge \tau_{2i+1} : S_t \in \partial D(0,R)\}.$$

Abbreviating $\bar{\tau} = T_{\partial D(0,R')}$, note that $\bar{\tau} = \tau_{2I+1}$ for some (unique) nonnegative integer $I$. For any $i \ge 1$, we can write $\{B, I = i\} = \{B_i, \tau_{2i} < \bar{\tau}\} \cap (\{I = 0\} \circ \theta_{\tau_{2i}})$ for some $B_i \in \mathcal{F}_{\tau_{2i}}$, so by the strong Markov property at $\tau_{2i}$,

$$\mathbb{E}^z[X_{\bar{\tau}} = w; B, I = i] = \mathbb{E}^z[\mathbb{E}^{X_{\tau_{2i}}}(X_{\bar{\tau}} = w, I = 0); B_i, \tau_{2i} < \bar{\tau}]$$

and

$$\mathbf{P}^z(B, I = i) = \mathbb{E}^z[\mathbb{E}^{X_{\tau_{2i}}}(I = 0); B_i, \tau_{2i} < \bar{\tau}].$$

Consequently, for all $i \geq 1$,

(2.17)
$$\mathbb{E}^z[X_{\bar{\tau}} = w; B, I = i]$$
$$\geq \mathbf{P}^z(B, I = i) \min_{x \in \partial D(0,R)} \frac{\mathbb{E}^x(X_{\bar{\tau}} = w; I = 0)}{\mathbb{E}^x(I = 0)}.$$

Necessarily $\mathbf{P}^z(B | I = 0) \in \{0, 1\}$ and is independent of $z$ for any $B \in \mathcal{H}$, implying that (2.17) applies for $i = 0$ as well. By our assumptions about $r, R, R'$, (2.4), (2.5) and (2.6) there exists $c < \infty$ such that for any $z, x \in \partial D(0,R)$ and $w \in \partial D(0,R')$,

$$\mathbb{E}^x(X_{\bar{\tau}} = w; I = 0) \geq (1 - cR/R')\mathbb{E}^x(I = 0)H_{\partial D(0,R')}(z, w).$$

Hence, summing (2.17) over $I = 0, 1, \ldots$, we get that

$$\mathbb{E}^z[X_{\bar{\tau}} = w, B] \geq (1 - cR/R')\mathbf{P}^z(B)H_{\partial D(0,R')}(z, w).$$

A similar argument shows that

$$\mathbb{E}^z[X_{\bar{\tau}} = w, B] \leq (1 + cR/R')\mathbf{P}^z(B)H_{\partial D(0,R')}(z, w),$$

and we thus obtain (2.15).

By the Markov property at $\tau_1$, for any $z \in \partial D(0, R)$,

$$\mathbf{P}^z(B) = \mathbf{P}^z(B, I = 0)$$
$$+ \sum_{x \in \partial D(0,r)} \widehat{H}_{\partial D(0,r) \cup \partial D(0,R')}(z, x)\mathbf{P}^x(B).$$

The term involving $\{B, I = 0\}$ is dealt with by (2.4), and (2.16) follows by (2.12) and our assumptions about $r, R$ and $R'$ values. $\quad\square$

Building upon Lemma 2.3 we quantify the independence between the $\sigma$-algebra $\mathcal{G}^x$ of excursions from $\partial D(x, R')$ to $\partial D(x, R)$, and the $\sigma$-algebra $\mathcal{H}^x(m)$ of excursions from $\partial D(x, r)$ to $\partial D(x, R)$ which occur during the first $m$ excursions from $\partial D(x, R)$ to $\partial D(x, R')$. To this end, fix $4r < R < R' < K/2$ and $x \in \mathbb{Z}_K^2$, let $\overline{\tau}_0 = 0$, and for $i = 1, 2, \ldots$ define

$$\tau_i = \inf\{t \geq \overline{\tau}_{i-1} : X_t \in \partial D(x, R)\},$$
$$\overline{\tau}_i = \inf\{t \geq \tau_i : X_t \in \partial D(x, R')\}.$$

Then $\mathcal{G}^x$ is the $\sigma$-algebra generated by the excursions $\{e^{(j)}, j = 1, \ldots\}$, where $e^{(j)} = \{X_t : \overline{\tau}_{j-1} \leq t \leq \tau_j\}$ is the $j$th excursion from $\partial D(x, R')$ to $\partial D(x, R)$ (so for $j = 1$ we do begin at $t = 0$). We denote by $\mathcal{H}^x(m)$ the $\sigma$-algebra generated by all excursions from $\partial D(x, r)$ to $\partial D(x, R)$ from time $\tau_1$ until time $\overline{\tau}_m$. In more detail, for each $j = 1, 2, \ldots, m$ let $\overline{\zeta}_{j,0} = \tau_j$ and for $i = 1, \ldots$ define

$$\zeta_{j,i} = \inf\{t \geq \overline{\zeta}_{j,i-1} : X_t \in \partial D(x, r)\},$$

$$\overline{\zeta}_{j,i} = \inf\{t \geq \zeta_{j,i} : X_t \in \partial D(x, R)\}.$$

Let $v_{j,i} = \{X_t : \zeta_{j,i} \leq t \leq \overline{\zeta}_{j,i}\}$ and $Z^j = \sup\{i \geq 0 : \overline{\zeta}_{j,i} < \overline{\tau}_j\}$. Then, $\mathcal{H}^x(m)$ is the product $\sigma$-algebra generated by the $\sigma$-algebras $\mathcal{H}^x_j = \sigma(v_{j,i}, i = 1, \ldots, Z^j)$ of the excursions between times $\tau_j$ and $\overline{\tau}_j$, for $j = 1, \ldots, m$.

LEMMA 2.4. *There exists $C < \infty$ such that uniformly over $\sqrt{R} < 4r <$* $R < R' < K/2$ *with $\log(R'/R) \geq (1/4)\log(R/r),$ all $m \leq R/(r\log(R/r)),$* $x, y_0, y_1 \in \mathbb{Z}_K^2$ *and $A \in \mathcal{H}^x(m),$*

$$
\begin{aligned}
\left(1 - Cm\frac{r}{R}\log\frac{R}{r}\right)\mathbf{P}^{y_1}(A) &\leq \mathbf{P}^{y_0}(A|\mathcal{G}^x) \\
(2.18) \qquad\qquad &\leq \left(1 + Cm\frac{r}{R}\log\frac{R}{r}\right)\mathbf{P}^{y_1}(A).
\end{aligned}
$$

PROOF. Applying the monotone class theorem to the algebra of their finite disjoint unions, it suffices to prove (2.18) for the generators of the product $\sigma$-algebra $\mathcal{H}^x(m)$ of the form $A = A_1 \times A_2 \times \cdots \times A_m$, with $A_j \in \mathcal{H}^x_j$ for $j = 1, \ldots, m$. Conditioned upon $\mathcal{G}^x$ the events $A_j$ are independent. Further, each $A_j$ then has the conditional law of an event $B_j$ in the $\sigma$-algebra $\mathcal{H}$ of Lemma 2.3, for some random $z_j = X_{\tau_j} - x \in \partial \tilde{D}(0, R)$ and $w_j = X_{\overline{\tau}_j} - x \in \partial D(0, R')$, both measurable on $\mathcal{G}^x$. By our conditions on $r$, $R$ and $R'$, the uniform estimates (2.15) and (2.16) yield that for any fixed $z' \in \partial D(0, R)$,

$$
\begin{aligned}
\mathbf{P}^{y_0}(A_1 \times A_2 \times \cdots \times A_m | \mathcal{G}^x) &= \prod_{j=1}^m \mathbf{P}^{z_j}(B_j | X_{T_{\partial D(0,R')}} = w_j) \\
(2.19) \qquad\qquad &= \prod_{j=1}^m \left(1 + O\left(\frac{R}{R'}\right)\right)\mathbf{P}^{z_j}(B_j) \\
&= \left(1 + O\left(\frac{r}{R}\log\frac{R}{r}\right)\right)^m \prod_{j=1}^m \mathbf{P}^{z'}(B_j).
\end{aligned}
$$

Since $m \leq R/(r\log(R/r))$ and the right-hand side of (2.19) depends neither on $y_0 \in \mathbb{Z}_K^2$ nor on the extra information in $\mathcal{G}^x$, we get (2.18) by averaging over $\mathcal{G}^x$. □

REMARK. Lemma 2.3, which deals with the path of the walk in $\overline{D}(0, R')$, applies for the simple random walk $S_n$ in $\mathbb{Z}^2$. Consequently, by the same argument as above, the bounds of (2.18) also apply for $S_n$.

**3. Hitting time estimates and upper bounds.** For any first hitting time $T$ we set $\|T\| = \sup_y \mathbb{E}^y(T)$. By Kac's moment formula for the strong Markov process $X_n$ (see [5], (6)), we have for any $n$ and $y$

$$\mathbb{E}^y(T^n) \leq n! \mathbb{E}^y(T) \|T\|^{n-1}. \tag{3.1}$$

Throughout this section, consider constants $r, R$ such that $0 < 2r < R \leq \frac{1}{2}K$. Per fixed $x \in \mathbb{Z}_K^2$, we let

$$\tau^{(0)} = \inf\{t \geq 0 : X_t \in \partial D(x, r)\}, \tag{3.2}$$

$$\sigma^{(1)} = \inf\{t \geq 0 : X_{t+\tau^{(0)}} \in \partial D(x, R)\}, \tag{3.3}$$

and define inductively for $j = 1, 2, \ldots$

$$\tau^{(j)} = \inf\{t \geq \sigma^{(j)} : X_{t+\mathfrak{T}_{j-1}} \in \partial D(x, r)\}, \tag{3.4}$$

$$\sigma^{(j+1)} = \inf\{t \geq 0 : X_{t+\mathfrak{T}_j} \in \partial D(x, R)\}, \tag{3.5}$$

where $\mathfrak{T}_j = \sum_{i=0}^j \tau^{(i)}$ for $j = 0, 1, 2, \ldots$. Thus $\tau^{(j)}$, $j \geq 1$, is the length of the $j$th excursion $\mathcal{E}_j$ from $\partial D(x, r)$ to itself via $\partial D(x, R)$, and $\sigma^{(j)}$ is the amount of time it takes to hit $\partial D(x, R)$ during the $j$th excursion $\mathcal{E}_j$. Hereafter, we set $\tau = \tau^{(1)}$ and use the abbreviation $\partial r = \partial D(x, r)$.

The following lemma will be used repeatedly.

LEMMA 3.1. *There exists $c_1 < \infty$ such that for all $1 \geq \eta \geq c_1(1/r + r/R)$ and $R \leq K/6$,*

$$(1 - \eta)\frac{2}{\pi}K^2 \log\left(\frac{R}{r}\right) \leq \min_{x, y \in \mathbb{Z}_K^2} \mathbb{E}^y(\tau) \leq \max_{x, y \in \mathbb{Z}_K^2} \mathbb{E}^y(\tau)$$

$$\leq (1 + \eta)\frac{2}{\pi}K^2 \log\left(\frac{R}{r}\right), \tag{3.6}$$

$$\max_{x \in \mathbb{Z}_K^2} \max_{y \in \partial D(x, R)} \mathbb{E}^y(T_{\partial D(x, r)}) \leq c_1 K^2 \log\left(\frac{R}{r}\right), \tag{3.7}$$

*and for all $r \geq c_1$,*

$$\max_{x \in \mathbb{Z}_K^2} \|T_{\partial D(x, r)}\| \leq c_1 K^2 \log\left(\frac{K}{r}\right). \tag{3.8}$$

PROOF. Let $X_0$ be distributed uniformly on $\mathbb{Z}_K^2$. Then $\{X_t\}$ is a stationary and ergodic stochastic process. By Birkhoff's ergodic theorem we then have that

$$\lim_{T \to \infty} \frac{1}{T} \sum_{i=0}^T \mathbf{1}_{\{x\}}(X_i) = \frac{1}{K^2} \qquad \text{a.s.}$$

Thus, with $\mathfrak{T}_{-1} = 0$,

$$(3.9) \qquad \lim_{n \to \infty} \frac{1/n \sum_{j=0}^{n} \sum_{i=0}^{\tau^{(j)}} \mathbf{1}_{\{x\}}(X_{i+\mathfrak{T}_{j-1}})}{1/n \sum_{j=0}^{n} \tau^{(j)}} = \frac{1}{K^2} \qquad \text{a.s.}$$

For $j \geq 1$ set $Z_j = \tau^{(j)} - \mathbb{E}^{\rho}(\tau^{(j)}|\mathcal{F}_{\mathfrak{T}_{j-1}}) = \tau^{(j)} - \mathbb{E}^{X_{\mathfrak{T}_{j-1}}}(\tau)$, where $\rho$ is uniform measure on $\mathbb{Z}_K^2$. By the strong Markov property we see that $\{Z_j\}$ is an orthogonal sequence. Since any irreducible Markov chain with finite state space is positive recurrent, we have that $\|T_{\partial r}\|, \|T_{\partial R}\| < \infty$, and using (3.1) we see that the sequence $\{\tau^{(j)}\}$ and hence $\{Z_j\}$ has uniformly bounded second moments. It follows from Rajchman's strong law of large numbers (see, e.g., [2], Theorem 5.1.2) that

$$(3.10) \qquad \lim_{n \to \infty} \frac{1}{n} \sum_{j=1}^{n} \{\tau^{(j)} - \mathbb{E}^{X_{\mathfrak{T}_{j-1}}}(\tau)\} = 0 \qquad \text{a.s.}$$

Similarly, set $\sigma^{(0)} = \tau^{(0)}$ and for $j \geq 0$ let

$$Y_j = \sum_{i=0}^{\tau^{(j)}} \mathbf{1}_{\{x\}}(X_{i+\mathfrak{T}_{j-1}}) = \sum_{i=0}^{\sigma^{(j)}} \mathbf{1}_{\{x\}}(X_{i+\mathfrak{T}_{j-1}}),$$

$$\widetilde{Y}_j = Y_j - \mathbb{E}^{\rho}(Y_j|\mathcal{F}_{\mathfrak{T}_{j-1}}) = Y_j - \mathbb{E}^{X_{\mathfrak{T}_{j-1}}}(Y_1).$$

By the strong Markov property $\{\widetilde{Y}_j\}$ is also an orthogonal sequence, and since $Y_j \leq \tau^{(j)}$, the sequence $\{\widetilde{Y}_j\}$ also has uniformly bounded second moments. Thus, by Rajchman's strong law of large numbers,

$$(3.11) \qquad \lim_{n \to \infty} \frac{1}{n} \sum_{j=1}^{n} \{Y_j - \mathbb{E}^{X_{\mathfrak{T}_{j-1}}}(Y_1)\} = 0 \qquad \text{a.s.}$$

It follows from (2.2) that for some finite universal constant $c_0 \geq 1$ and all $1 \leq r \leq R/3 \leq K/6$,

$$(3.12) \qquad \begin{aligned} \frac{2}{\pi} \log\left(\frac{R}{r}\right) - c_0 r^{-1} &\leq \min_x \min_{y \in \partial r} \mathbb{E}^y(Y_1) \\ &\leq \max_x \max_{y \in \partial r} \mathbb{E}^y(Y_1) \\ &\leq \frac{2}{\pi} \log\left(\frac{R}{r}\right) + c_0 r^{-1}. \end{aligned}$$

With $\tau^{(0)}$ finite, we get by combining (3.9), (3.10) and (3.11) that almost surely,

$$\lim_{n \to \infty} \frac{(1/n) \sum_{j=1}^{n} \mathbb{E}^{X_{\mathfrak{T}_{j-1}}}(\tau)}{(1/n) \sum_{j=1}^{n} \mathbb{E}^{X_{\mathfrak{T}_{j-1}}}(Y_1)} = K^2.$$

Consequently, in view of (3.12), for some finite universal constant $c_1$ and all $1 > \eta \geq c_1(1/r + r/R)$,

$$(3.13) \quad \frac{2}{\pi}\left(1 - \frac{\eta}{3}\right)K^2 \log\left(\frac{R}{r}\right) \leq \max_{y \in \partial r} \mathbb{E}^y(\tau),$$

$$\min_{y \in \partial r} \mathbb{E}^y(\tau) \leq \frac{2}{\pi}\left(1 + \frac{\eta}{3}\right)K^2 \log\left(\frac{R}{r}\right).$$

For $y \in \partial r$, we have $\tau^{(0)} = 0$ and by the strong Markov property at the stopping time $\sigma^{(1)}$,

$$(3.14) \quad \mathbb{E}^y(\tau) = \mathbb{E}^y(T_{\partial R}) + \sum_{z \in \partial R} H_{\partial R}(y, z) \mathbb{E}^z(T_{\partial r}).$$

Thus, enlarging $c_0$ as needed, it follows from (2.3) and Lemma 2.1 that for all $1 \leq r \leq R/c_0$,

$$(3.15) \quad \max_{y \in \partial r} \mathbb{E}^y(\tau) \leq \left(1 + c_0 \frac{r}{R}\right) \min_{y \in \partial r} \mathbb{E}^y(\tau).$$

Taking also $c_1 \geq 3c_0$, we get (3.6) by combining (3.13) and (3.15).

Turning to prove (3.7), consider (3.14) for $y \in \partial r$ and $3R$ instead of $R$. Then, by (3.6) and (2.9),

$$(3.16) \quad c(1/3)\mathbb{E}^{\mu_{3R}}(T_{\partial r}) \leq 2K^2 \log(3R/r).$$

Using the strong Markov property, (2.3), (2.9) and (3.16), we thus have that for any $y \in \partial R$,

$$\begin{aligned}(3.17) \quad \mathbb{E}^y(T_{\partial r}) &\leq \mathbb{E}^y(T_{\partial 3R}) + \mathbb{E}^y(T_{\partial r} - T_{\partial 3R}; T_{\partial r} > T_{\partial 3R}) \\ &\leq (3R + 1)^2 + C\left(\frac{1}{3}\right)\mathbb{E}^{\mu_{3R}}(T_{\partial r}) \\ &\leq c_2 K^2 \log\left(\frac{R}{r}\right),\end{aligned}$$

for some universal $c_2 < \infty$ and any $r$, $R$ as in the statement of (3.7). Making sure that $c_1 \geq c_2$, this completes the proof of (3.7).

To prove (3.8) we use the bound (3.17) when the distance of $y$ from $x$ is between $R_0 = r/c_1$ and $K/6$, and that of (2.3) when $y \in D(x, r)$. As for $y \in D(x, R_0) \setminus D(x, r)$, since

$$\mathbb{E}^y(T_{\partial r}) \leq \mathbb{E}^y(T_{\partial R_0}) + \max_{z \in \partial R_0} \mathbb{E}^z(T_{\partial r}),$$

we get the stated bound by combining (2.3) (for the first term above) and (3.17). Finally, fixing $y \in \mathbb{Z}_K^2 \setminus D(x, K/6)$, we establish the bound of (3.8) by noticing that the value of $\mathbb{E}^y(T_{\partial r})$ for the random walk on $\mathbb{Z}_{2^\ell K}^2$ is then nondecreasing in $\ell$, and adjusting $c_1$ accordingly (to accommodate the use of, say, $\mathbb{Z}_{16K}^2$). $\quad\square$

The following lemma, which shows that excursion times are concentrated around their mean, will be used to relate excursions to hitting times.

LEMMA 3.2. *With the above notation, we can find $\delta_0 > 0$ and $C > 0$ such that if $R \le K/2$ and $\delta \le \delta_0$ with $\delta \ge 6c_1(1/r + r/R)$, then for all $x, x_0 \in \mathbb{Z}_K^2$,*

$$(3.18) \quad \mathbf{P}^{x_0}\left(\sum_{j=0}^{N} \tau^{(j)} \le (1-\delta)\frac{2K^2 \log(R/r)}{\pi} N\right)$$

$$\le e^{-C\delta^2(\log(R/r)/\log(K/r))N}$$

*and*

$$(3.19) \quad \mathbf{P}^{x_0}\left(\sum_{j=0}^{N} \tau^{(j)} \ge (1+\delta)\frac{2K^2 \log(R/r)}{\pi} N\right)$$

$$\le e^{-C\delta^2(\log(R/r)/\log(K/r))N}.$$

PROOF. With $\tau = \tau^{(1)} = \{T_{\partial R} + T_{\partial r} \circ \theta_{T_{\partial R}}\} \circ \theta_{T_{\partial r}}$, clearly

$$\max_y \mathbb{E}^y(\tau^n) \le \max_{y \in \partial r} \mathbb{E}^y(\{T_{\partial R} + T_{\partial r} \circ \theta_{T_{\partial R}}\}^n)$$

$$\le \sum_{j=0}^n \binom{n}{j} \max_{y \in \partial r} \mathbb{E}^y(T_{\partial R}^j (T_{\partial r}^{n-j} \circ \theta_{T_{\partial R}}))$$

$$\le \sum_{j=0}^n \binom{n}{j} \max_{y \in \partial r} \mathbb{E}^y(T_{\partial R}^j) \max_{z \in \partial R} \mathbb{E}^z(T_{\partial r}^{n-j}).$$

Let $v = \frac{2K^2}{\pi} \log(R/r)$ and $u = \frac{2K^2}{\pi} \log(K/r)$. Thus, by (3.1) and (3.7), there exists a universal constant $c_3 < \infty$ such that for all $x \in \mathbb{Z}_K^2$,

$$\max_y \mathbb{E}^y(\tau^n) \le \max_{y \in \partial r} \mathbb{E}^y(T_{\partial R}) \|T_{\partial R}\|^{n-1}$$

$$(3.20) \qquad\qquad + 2c_1 \sum_{j=0}^{n-1} n! \|T_{\partial R}\|^j v \|T_{\partial r}\|^{n-j-1}$$

$$\le v(c_3 u)^{n-1}(n+1)!,$$

where we also used (2.3) and (3.8) in the last inequality. Taking $\eta = \delta/6 > 0$, with our choice of $r$ and $R$, it thus follows by (3.6) that for $\rho = c_4 uv$ and all

$\theta \geq 0$,

$$
\begin{aligned}
\max_x \max_{y \in \partial D(x,r)} \mathbb{E}^y(e^{-\theta \tau}) &\leq 1 - \theta \min_x \min_{y \in \partial D(x,r)} \mathbb{E}^y(\tau) \\
(3.21) \qquad\qquad &+ \frac{\theta^2}{2} \max_x \max_{y \in \partial D(x,r)} \mathbb{E}^y(\tau^2) \\
&\leq 1 - \theta(1-\eta)v + \rho\theta^2 \\
&\leq \exp(\rho\theta^2 - \theta(1-\eta)v).
\end{aligned}
$$

Since $\tau^{(0)} \geq 0$, using Chebyshev's inequality we bound the left-hand side of (3.18) by

$$
\begin{aligned}
(3.22) \qquad \mathbf{P}^{x_0}\left(\sum_{j=1}^N \tau^{(j)} \leq (1-6\eta)vN\right) &\leq e^{\theta(1-3\eta)vN}\mathbb{E}^{x_0}(e^{-\theta\sum_{j=1}^N \tau^{(j)}}) \\
&\leq e^{-\theta vN\delta/3}\left[e^{\theta(1-\eta)v}\max_{y \in \partial D(x,r)} \mathbb{E}^y(e^{-\theta\tau})\right]^N,
\end{aligned}
$$

where the last inequality follows by the strong Markov property of $X_t$ at $\{\mathfrak{T}_j\}$. Combining (3.21) and (3.22) for $\theta = \delta v/(6\rho)$, results in (3.18) with $C = 1/(36c_4)$.

Since $\tau^{(0)} = T_{\partial r}$, by (3.1) and (3.8) there exist universal constants $c_5, c_6 < \infty$ such that

$$
\max_{x,y} \mathbb{E}^y(e^{\tau^{(0)}/c_5 u}) \leq c_6,
$$

implying that

$$
\mathbf{P}^{x_0}\left(\tau^{(0)} \geq \frac{\delta}{3}vN\right) = \mathbf{P}^{x_0}\left(\frac{\tau^{(0)}}{c_5 u} \geq \frac{\delta}{3c_5}\frac{v}{u}N\right) \leq c_6 e^{-(3c_5)^{-1}\delta(v/u)N}.
$$

Thus, the proof of (3.19), in analogy with that of (3.18), comes down to bounding

$$
\mathbf{P}^{x_0}\left(\sum_{j=1}^N \tau^{(j)} \geq (1+4\eta)vN\right) \leq e^{-\theta\delta vN/3}\left(e^{-\theta(1+2\eta)v}\max_{y \in \partial D(x,r)} \mathbb{E}^y(e^{\theta\tau})\right)^N.
$$

Noting that, by (3.20) and (3.6), there exists a universal constant $c_8 < \infty$ such that for $\rho = c_8 uv$ and all $0 < \theta < 1/(2c_3 u)$,

$$
\begin{aligned}
\max_x \max_{y \in \partial D(x,r)} \mathbb{E}^y(e^{\theta\tau}) &\leq 1 + \theta(1+\eta)v + \max_x \max_{y \in \partial D(x,r)} \sum_{n=2}^\infty \frac{\theta^n}{n!}\mathbb{E}^y(\tau^n) \\
&\leq 1 + \theta(1+2\eta)v + \rho\theta^2 \leq \exp(\theta(1+2\eta)v + \rho\theta^2).
\end{aligned}
$$

Taking $\delta_0 < 3c_8/c_3$, the proof of (3.19) now follows that of (3.18). $\quad\square$

We next apply Lemma 3.2 to bound the upper tail of $\mathcal{T}_K(x)$, the first hitting time of $x \in \mathbb{Z}_K^2$.

LEMMA 3.3. *For any $\delta > 0$ we can find $c < \infty$ and $K_0 < \infty$ so that for all $K \geq K_0$, $y \geq 0$ and $x, x_0 \in \mathbb{Z}_K^2$,*

$$(3.23) \qquad \mathbf{P}^{x_0}(\mathcal{T}_K(x) \geq y(K \log K)^2) \leq cK^{-(1-\delta)\pi y/2}.$$

PROOF. Fix $\delta \in (0, \delta_0)$. Set $R = K/7$ and $r = R/\log K$, noting that Lemma 3.2 then applies for all $K \geq K_0$ and some $K_0 = K_0(\delta) < \infty$. Fixing $y \geq 0$ and such $K$, let

$$n_K := (1-\delta)\frac{\pi y (\log K)^2}{2 \log(R/r)}.$$

Then,

$$
\begin{aligned}
(3.24) \qquad & \mathbf{P}^{x_0}(\mathcal{T}_K(x) \geq y(K \log K)^2) \\
& \leq \mathbf{P}^{x_0}\left(\mathcal{T}_K(x) \geq \sum_{j=0}^{n_K} \tau^{(j)}\right) + \mathbf{P}^{x_0}\left(\sum_{j=0}^{n_K} \tau^{(j)} \geq y(K \log K)^2\right).
\end{aligned}
$$

It follows from (3.19) that

$$\mathbf{P}^{x_0}\left(\sum_{j=0}^{n_K} \tau^{(j)} \geq y(K \log K)^2\right) \leq e^{-C'y(\log K)^2/\log \log K}$$

for some $C' = C'(\delta) > 0$. Moreover, the first probability in (3.24) is bounded above by the probability of not hitting $x$ during $n_K$ excursions of SRW in $\mathbb{Z}^2$, each starting at some point in $\partial D(x, r)$ and ending at $\partial D(x, R)$, so that by (2.1)

$$
\begin{aligned}
(3.25) \qquad & \mathbf{P}^{x_0}\left(\mathcal{T}_K(x) \geq \sum_{j=0}^{n_K} \tau^{(j)}\right) \leq \left(1 - \frac{\log R/r + O(1/\log K)}{\log R}\right)^{n_K} \\
& \leq e^{-(1-2\delta)\log(K)\pi y/2},
\end{aligned}
$$

and (3.23) follows. $\square$

We next provide the required upper bounds in Proposition 1.1. Namely, for any $\alpha \in (0, 1]$ and $\gamma > 0$, we have by Lemma 3.3, that for $\gamma/(2\alpha) > \delta > 0$ small enough,

$$
\begin{aligned}
(3.26) \qquad & \mathbf{P}\left(\left|\left\{x \in \mathbb{Z}_K^2 : \frac{\mathcal{T}_K(x)}{(K \log K)^2} \geq 4\alpha/\pi\right\}\right| \geq K^{2(1-\alpha)+\gamma}\right) \\
& \leq K^{-2(1-\alpha)-\gamma}\mathbb{E}\left(\left|\left\{x \in \mathbb{Z}_K^2 : \frac{\mathcal{T}_K(x)}{(K \log K)^2} \geq 4\alpha/\pi\right\}\right|\right)
\end{aligned}
$$

$$= K^{-2(1-\alpha)-\gamma} \sum_{x \in \mathbb{Z}_K^2} \mathbf{P}\left(\frac{\mathcal{T}_K(x)}{(K \log K)^2} \geq 4\alpha/\pi\right)$$

$$\leq K^{2\delta\alpha - \gamma} \underset{K \to \infty}{\longrightarrow} 0.$$

**4. Lower bounds for probabilities.** Fixing $a < 2$, we prove in this section that for any $\delta > 0$ there exists $n_0(\delta) < \infty$ such that

$$(4.1) \qquad \mathbf{P}\left(\left|\left\{x \in \mathbb{Z}_{K_n}^2 : \frac{\mathcal{T}_{K_n}(x)}{(K_n \log K_n)^2} \geq 2a/\pi\right\}\right| \geq K_n^{2-a-\delta}\right) \geq 1 - 2\delta,$$

for all integers $K_n = n^{\bar{\gamma}}(n!)^3$ with $n \geq n_0$ and $\bar{\gamma} \in \mathcal{I} = [b, b+4]$ for some universal $b \geq 10$ (determined in Lemma 4.2). Because such $K_n$ cover all large enough integers, it follows from (4.1) that

$$\lim_{m \to \infty} \mathbf{P}\left(\left|\left\{x \in \mathbb{Z}_m^2 : \frac{\mathcal{T}_m(x)}{(m \log m)^2} \geq 2a/\pi\right\}\right| \geq m^{2-a-\delta}\right) = 1,$$

which in view of (3.26) results with Proposition 1.1. Hereafter, any estimate involving the fixed sequence $K_n = n^{\bar{\gamma}}(n!)^3$ holds uniformly in $\bar{\gamma} \in \mathcal{I}$ (even if this is not stated explicitly). Consequently, we may and shall prove each of our results only for this sequence, which already implies that they hold true for all integers large enough.

We start by constructing a subset of the set appearing in (4.1), the probability of which is easier to bound below. To this end, let $r_0 = 0$ and $r_k = (k!)^3, k = 1, \ldots$ For any $a > 0$ set $n_k = n_k(a) = 3ak^2 \log k$ and for $x \in \mathbb{Z}_{K_n}^2$ and $k = 3, \ldots, n$, let $\mathcal{R}_k^x = \mathcal{R}_k^x(a)$ denote the time until completion of the first $n_k(a)$ excursions from $\partial D(x, r_{k-1})$ to $\partial D(x, r_k)$. (In the notation of Section 3, if we set $R = r_k$ and $r = r_{k-1}$, then $\mathcal{R}_k^x = \sum_{j=0}^{n_k} \tau^{(j)}$.) For $x \in \mathbb{Z}_{K_n}^2$, $2 \leq l \leq k-1$, let $N_{k,l}^x = N_{k,l}^x(a)$ denote the number of excursions from $\partial D(x, r_{l-1})$ to $\partial D(x, r_l)$ until time $\mathcal{R}_k^x(a)$. Let $N_{k,0}^x = N_{k,0}^x(a)$ denote the number of visits to $x$ prior to time $\mathcal{R}_k^x(a)$.

Fix $\rho < (2-a)/2$. Writing $m \overset{k}{\sim} n_k$ if $|m - n_k| \leq k$, we will say that a point $x \in \mathbb{Z}_{K_n}^2$ is *n-successful* if

$$(4.2) \qquad N_{n,0}^x = 0, \qquad N_{n,k}^x \overset{k}{\sim} n_k \qquad \forall k = \rho n, \ldots, n-1.$$

Note that $N_{n,0}^x = 0$ is equivalent to the statement $\mathcal{T}_{K_n}(x) > \mathcal{R}_n^x$. Hence the next lemma relates the notions of $n$-successful and first hitting times.

LEMMA 4.1. *Let*

$$\mathcal{S}_n = \{x \in \mathbb{Z}_{K_n}^2 : \mathcal{T}_{K_n}(x) > \mathcal{R}_n^x\}.$$

*Then, for some $c > 0$ independent of $\bar\gamma$ and all $n \geq n_0$,*

$$\mathbf{P}\left( \bigcup_{x \in \mathcal{S}_n} \left\{ \frac{\mathcal{T}_{K_n}(x)}{(K_n \log K_n)^2} \leq 2a/\pi - 2/\log n \right\} \right) \leq c^{-1} e^{-cn^2/\log n}.$$

PROOF. We have that for some $C > 0$ and $n_0 < \infty$, both independent of $\bar\gamma$, all $n \geq n_0$ and any $x, x_0 \in \mathbb{Z}^2_{K_n}$,

$$P_x := \mathbf{P}^{x_0}(\mathcal{T}_{K_n}(x) \leq (2a/\pi - 2/\log n)(K_n \log K_n)^2, \mathcal{T}_{K_n}(x) > \mathcal{R}^x_n)$$

$$\leq \mathbf{P}^{x_0}\left( \sum_{j=0}^{3an^2 \log n} \tau^{(j)} \leq (2a/\pi - 1/\log n) K_n^2 (3n \log n)^2 \right) \leq e^{-Cn^2/\log n},$$

where the last inequality is an application of (3.18) with $R = r_n$, $r = r_{n-1}$ [so $\log(R/r) = 3 \log n$] and $\delta = \pi/(2a \log n)$. To complete the proof of the lemma, sum over $x \in \mathbb{Z}^2_{K_n}$ and let $c < C/2$ be such that $c^{-1} e^{-cn_0^2} \geq 1$. □

For any $x \in \mathbb{Z}^2_{K_n}$ let $Y(n, x)$ be the indicator random variable for the event $\{x$ is $n$-successful$\}$. In view of Lemma 4.1, we have (4.1) as soon as we show that

$$(4.3) \qquad \mathbf{P}\left( \sum_{x \in \mathbb{Z}^2_{K_n}} Y(n, x) \geq K_n^{2-a-\delta} \right) \geq 1 - \delta,$$

for any $\delta > 0$, all $n$ sufficiently large and $\bar\gamma \in \mathcal{I}$.

Adopting hereafter the convention that $o(1_n)$ terms are uniform in $\bar\gamma \in \mathcal{I}$, the key to the proof of (4.3) is the next lemma (whose proof is deferred to Section 5).

LEMMA 4.2. *Fix $\rho < \rho' < (2 - a)/2$ and let $l(x, y) = \max\{k : D(x, r_k + 1) \cap D(y, r_k + 1) = \varnothing\} \wedge n$. There exist $b \geq 10$ independent of $a$ and $\rho$, and $\bar{q}_n \geq r_n^{-a + o(1_n)}$ such that*

$$(4.4) \qquad \mathbf{P}(x \text{ is } n\text{-successful}) = (1 + o(1_n))\bar{q}_n,$$

*uniformly in $\bar\gamma \in \mathcal{I}$ and $x \in S_{K_n} := \mathbb{Z}^2_{K_n} \setminus D(0, r_n)$. Furthermore, for any $\varepsilon > 0$ we can find $C = C(b, \varepsilon) < \infty$ such that for all $n$ and any $x, y \in S_{K_n}$ with $\rho' n \leq l(x, y) < n$,*

$$(4.5) \qquad \mathbb{E}(Y(n, x)Y(n, y)) \leq \bar{q}_n^2 n^b C^{n - l(x, y)} \left( \frac{r_n}{r_{l(x, y)}} \right)^{a + \varepsilon},$$

*while for all $n$ and $x, y \in S_{K_n}$ with $l(x, y) = n$,*

$$(4.6) \qquad \mathbb{E}(Y(n, x)Y(n, y)) \leq (1 + o(1_n))\bar{q}_n^2.$$

Let

$$V_\ell = \sum_{x,y \in S_{K_n}, l(x,y)=\ell} \mathbb{E}(Y(n,x)Y(n,y)), \qquad \ell = 0, 1, \ldots, n.$$

Since, by (4.4),

$$\mathbb{E}\left(\sum_{x \in S_{K_n}} Y(n,x)\right) = (1+o(1_n))K_n^2 \bar{q}_n \geq K_n^{2-a+o(1_n)},$$

by (4.6) and the Paley–Zygmund inequality (see [6], page 8), inequality (4.3) is a direct consequence of the bound

$$(4.7) \qquad \sum_{\ell=0}^{n-1} V_\ell \leq o(1_n)K_n^4 \bar{q}_n^2.$$

Turning to prove (4.7), the definition of $l(x,y)$ implies that

$$d(x,y) < 2(r_{l(x,y)+1}+1),$$

and there are on $\mathbb{Z}_{K_n}^2$ at most $C_0 r_{\ell+1}^2$ points $y$ in the disc of radius $2(r_{\ell+1}+1)$ centered at $x$, where in the sequel we let $C_m$ denote generic finite constants that are independent of $n$. Since $2\rho' < 2-a$,

$$(4.8) \qquad \sum_{\ell=0}^{\rho'n-1} V_\ell \leq \sum_{x,y \in \mathbb{Z}_{K_n}^2, d(x,y) \leq 2r_{\rho'n}} \mathbb{E}(Y(n,x)) \leq C_1 \bar{q}_n K_n^2 r_{\rho'n}^2 \leq o(1_n)K_n^4 \bar{q}_n^2.$$

Choose $\varepsilon > 0$ such that $2-a-\varepsilon > 0$ and fix $\ell \in [\rho'n, n)$. Then, by (4.5), we have that

$$V_\ell \leq C_2 K_n^2 r_{\ell+1}^2 \bar{q}_n^2 n^b C^{n-\ell}\left(\frac{r_n}{r_\ell}\right)^{a+\varepsilon}.$$

Consequently,

$$\sum_{\ell=\rho'n}^{n-1} V_\ell \leq C_2 K_n^2 \bar{q}_n^2 n^b \sum_{\ell=\rho'n}^{n-1} C^{n-\ell} r_{\ell+1}^2 \left(\frac{r_n}{r_\ell}\right)^{a+\varepsilon}$$

$$(4.9) \qquad \leq C_2 \bar{q}_n^2 K_n^4 n^{-2\bar{\gamma}} n^{b+6} \sum_{\ell=\rho'n}^{n-1} C^{n-\ell}\left(\frac{r_\ell}{r_n}\right)^{2-a-\varepsilon}$$

$$\leq C_2 \bar{q}_n^2 K_n^4 n^{-2} \sum_{j=1}^\infty C^j r_j^{-(2-a-\varepsilon)}.$$

Combining (4.8) and (4.9) we establish (4.7), and hence complete the proof of (4.3) and thus of (4.1).

**5. First and second moment estimates.** For $y \in \mathbb{Z}_{K_n}^2$ and $n \geq l \geq 3$ let $\mathcal{G}_l^y$ denote the $\sigma$-algebra generated by the excursions of the random walk from $\partial D(y, r_l)$ to $\partial D(y, r_{l-1})$ as defined in Lemma 2.4 (for $R' = r_l$ and $R = r_{l-1}$). We start with the following corollary of Lemma 2.4 which plays a crucial role in the proof of Lemma 4.2.

COROLLARY 5.1. *Let* $\Gamma_l = \{N_{n,k}^y = m_k; k = 0, 2, \ldots, l-1\}$. *Then, uniformly over all* $n \geq l \geq n_0$, $\bar{\gamma} \in \mathcal{I}$, $m_l \overset{l}{\sim} n_l$, $\{m_k : k = 0, 2, \ldots, l-1\}$, $y \in \mathbb{Z}_{K_n}^2$ *and* $x_0, x_1 \in \mathbb{Z}_{K_n}^2 \setminus D(y, r_l)$,

$$
\mathbf{P}^{x_0}(\Gamma_l | N_{n,l}^y = m_l, \mathcal{G}_l^y) = (1 + O(l^{-1}(\log l)^2))\mathbf{P}^{x_1}(\Gamma_l | N_{n,l}^y = m_l)\mathbf{1}_{\{N_{n,l}^y = m_l\}}.
$$
(5.1)

PROOF. For $j = 1, 2, \ldots$ and $k = 2, \ldots, l-1$, let $Z_k^j$ denote the number of excursions from $\partial D(y, r_{k-1})$ to $\partial D(y, r_k)$ by the random walk during the time interval $[\tau_j, \overline{\tau}_j]$. Similarly, let $Z_0^j$ denote the number of visits to $y$ during this time interval. Clearly, the event

$$
A = \left\{ \sum_{j=1}^{m_l} Z_k^j = m_k : k = 0, 2, \ldots, l-1 \right\}
$$

belongs to the $\sigma$-algebra $\mathcal{H}^y(m_l)$ corresponding to $r = r_{l-2}$ in Lemma 2.4. It is easy to verify that starting at any $x \notin D(y, r_l)$, when the event $\{N_{n,l}^y = m_l\} \in \mathcal{G}_l^y$ occurs, it implies that $N_{n,k}^y = \sum_{j=1}^{m_l} Z_k^j$ for $k = 0, 2, \ldots, l-1$. Thus,

$$
\mathbf{P}^{x_0}(\Gamma_l | \mathcal{G}_l^y)\mathbf{1}_{\{N_{n,l}^y = m_l\}} = \mathbf{P}^{x_0}(A | \mathcal{G}_l^y)\mathbf{1}_{\{N_{n,l}^y = m_l\}}.
$$
(5.2)

For some universal constant $n_0 < \infty$ and all $l \geq n_0$ the conditions of Lemma 2.4 apply for our choice of $R' = r_l$, $R = r_{l-1}$ and $r = r_{l-2}$ with $(r/R)\log(R/r) \leq 4l^{-3}\log l$. With $m_l/(l^2 \log l)$ bounded above, by (2.18) we have, uniformly in $y \in \mathbb{Z}_{K_n}^2$ and $x_0, x_1 \in \mathbb{Z}_{K_n}^2 \setminus D(y, r_l)$,

$$
\mathbf{P}^{x_0}(A | \mathcal{G}_l^y) = (1 + O(l^{-1}(\log l)^2))\mathbf{P}^{x_1}(A).
$$
(5.3)

Hence,

$$
\mathbf{P}^{x_0}(\Gamma_l | \mathcal{G}_l^y)\mathbf{1}_{\{N_{n,l}^y = m_l\}} = (1 + O(l^{-1}(\log l)^2))\mathbf{P}^{x_1}(A)\mathbf{1}_{\{N_{n,l}^y = m_l\}}.
$$

Taking $x_0 = x_1$ and averaging, one has

$$
\mathbf{P}^{x_1}(\Gamma_l | N_{n,l}^y = m_l) = (1 + O(l^{-1}(\log l)^2))\mathbf{P}^{x_1}(A)
$$
(5.4)
$$
= (1 + O(l^{-1}(\log l)^2))\mathbf{P}^{x_0}(A | \mathcal{G}_l^y),
$$

where the second equality is due to (5.3). Using that $\{N_{n,l}^y = m_l\} \subset \mathcal{G}_l^y$, (5.2) and (5.4) imply (5.1). □

Proof of Lemma 4.2. We start by proving the first moment estimate (4.4). To this end, let $\bar{m} = (m_{\rho n}, m_{\rho n+1}, \ldots, m_n)$ be a candidate value of $N_{n,k}^x$, $k = \rho n, \ldots, n$, and set $|\bar{m}| = 2 \sum_{j=\rho n}^n m_j - 1$. Let $\mathcal{H}_n(\bar{m})$ be the collection of maps ("histories"),

$$s : \{1, 2, \ldots, |\bar{m}|\} \mapsto \{\rho n - 1, \rho n, \ldots, n\}$$

such that $s(1) = n - 1, s(|\bar{m}|) = n, |s(j+1) - s(j)| = 1$ and the number of up-crossings from $\ell - 1$ to $\ell$

$$u(\ell) =: |\{(j, j+1) | (s(j), s(j+1)) = (\ell-1, \ell)\}| = m_\ell.$$

The number of ways to partition the $u(\ell)$ up-crossings from $\ell - 1$ to $\ell$ before and among the $u(\ell+1)$ up-crossings from $\ell$ to $\ell + 1$ is

$$\binom{u(\ell+1) + u(\ell) - 1}{u(\ell)}.$$

Since the mapping $s$ is in one-to-one correspondence with the relative order of all its up-crossings,

$$|\mathcal{H}_n(\bar{m})| = \prod_{\ell = \rho n}^{n-1} \binom{m_{\ell+1} + m_\ell - 1}{m_\ell}.$$

To each path $\omega$ of the random walk $X$. we assign a "history" $h(\omega)$ as follows. Let $\tau(1)$ be the time of the first visit to $\partial D(x, r_{n-1})$, and define $\tau(2), \tau(3), \ldots$ to be the successive hitting times of different elements of $\{\partial D(x, r_{\rho n-1}), \ldots, \partial D(x, r_n)\}$. If $y \in \partial D(x, r_k)$ for some $k$, let $\Phi(y) = k$ and set $h(\omega)(j) = \Phi(\omega(\tau(j)))$. See Figure 1.

Let $h_{|_k}$ be the first $k$ coordinates of the sequence $h$. Let $p_\ell = \log(r_{\ell+1}/r_\ell)/\log(r_{\ell+1}/r_{\ell-1})$ and $q_\ell = \log(r_\ell/r_{\ell-1})/\log r_\ell$. Note that $\log(d(y, x)/r) = 1 + O(r^{-1})$ for any $r$, uniformly in $x$ and $y \in \partial D(x, r)$. So, applying the Markov property successively at the times $\tau(1), \tau(2), \ldots, \tau(|\bar{m}| - 1)$ and relying on (2.4) except for up-crossings from $\rho n - 1$ to $\rho n$, for which (2.1) applies, or for down-crossings from $n$ to $n - 1$, which occur with probability 1, we get that uniformly for any $s \in \mathcal{H}_n(\bar{m})$ and $x \in S_{K_n}$,

$$\mathbf{P}\{h_{|_{|\bar{m}|}} = s, \mathcal{T}_{K_n}(x) > \tau(|\bar{m}|)\}$$

$$= \prod_{\ell = \rho n}^{n-1} \{p_\ell + O(r_{\ell-1}^{-1})\}^{m_\ell} \{1 - p_\ell + O(r_{\ell-1}^{-1})\}^{m_{\ell+1}}$$

$$\times \{1 - q_{\rho n} + O((n \log n)^{-2})\}^{m_{\rho n}}.$$

Taking $m_n = n_n$, we see that uniformly in $x \in S_{K_n}$ and $\bar{\gamma} \in \mathcal{I}$,

$$
\begin{aligned}
\mathbf{P}(x \text{ is } n\text{-successful}) &= \sum_{\substack{m_{\rho n}, \ldots, m_{n-1} \\ |m_\ell - n_\ell| \le \ell}} \mathbf{P}\{h_{|_{|\bar{m}|}} \in \mathcal{H}_n(\bar{m}), \mathcal{T}_{K_n}(x) > \tau(|\bar{m}|)\} \\
&= (1 + o(1_n))\bar{q}_n,
\end{aligned}
$$

(5.5)

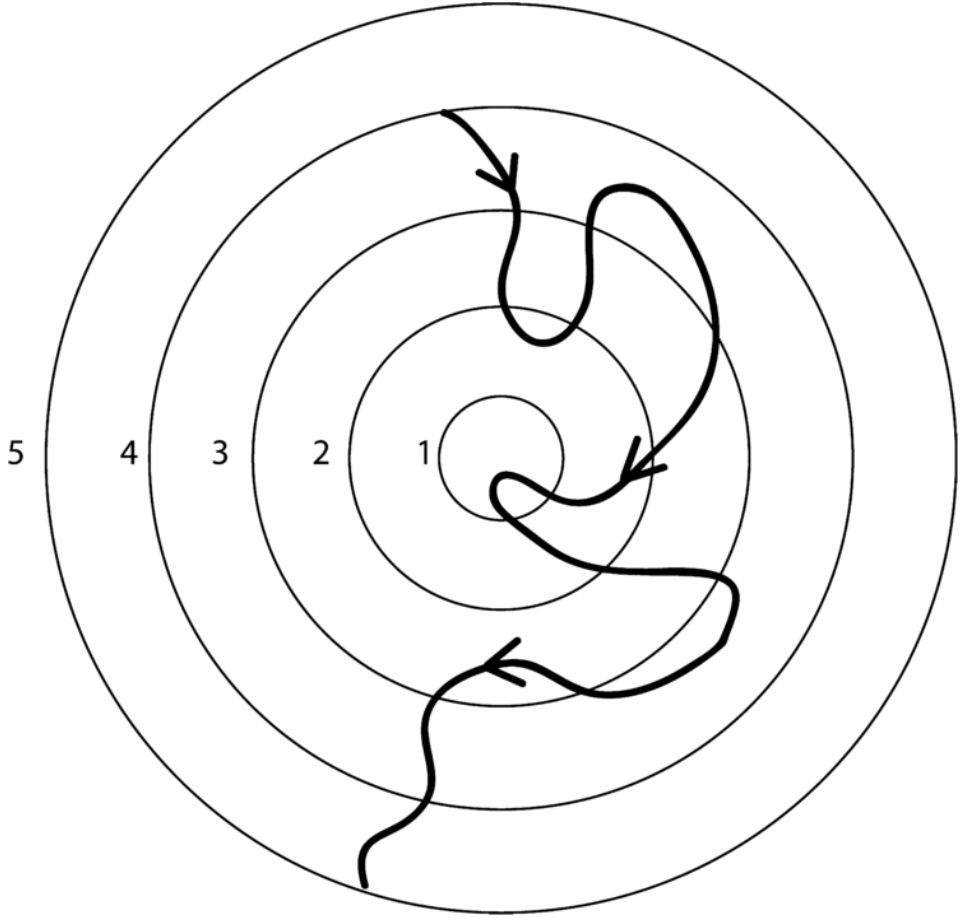

Fig. 1. *A path with "history" $h(\omega) = (4, 3, 2, 3, 2, 1, 2, 3, 4, 5)$.*

which is (4.4) for

$$(5.6) \quad \bar{q}_n = \sum_{\substack{m_{\rho n}, \dots, m_{n-1} \\ |m_\ell - n_\ell| \le \ell}} (1 - q_{\rho n})^{m_{\rho n}} \prod_{\ell = \rho n}^{n-1} \binom{m_{\ell+1} + m_\ell - 1}{m_\ell} p_\ell^{m_\ell} (1 - p_\ell)^{m_{\ell+1}}.$$

Since $p_\ell = 1/2 - O((\ell \log \ell)^{-1})$, by the proof of [3], Lemma 7.2, we have that uniformly in $m_\ell \overset{\ell}{\sim} n_\ell, m_{\ell+1} \overset{\ell+1}{\sim} n_{\ell+1}$

$$(5.7) \quad \begin{aligned} \frac{C' \ell^{-3a-1}}{\sqrt{\log \ell}} &\le \binom{m_{\ell+1} + m_\ell - 1}{m_\ell} p_\ell^{m_\ell} (1 - p_\ell)^{m_{\ell+1}} \\ &\le \frac{C \ell^{-3a-1}}{\sqrt{\log \ell}} \end{aligned}$$

with $0 < C', C < \infty$ independent of $\ell$. Further, with $q_\ell = \ell^{-1} + O(1/\ell \log \ell)$ we have that uniformly in $m_{\rho n} \overset{\rho n}{\sim} n_{\rho n}$

$$(5.8) \qquad (1 - q_{\rho n})^{m_{\rho n}} = r_{\rho n}^{-a + o(1_n)}.$$

Putting (5.6)–(5.8) together we see that $\bar{q}_n = r_n^{-a + o(1_n)}$, with the $o(1_n)$ term independent of $\bar{\gamma}$, as claimed.

Setting $M_l := \{l, l+1, \ldots, n-1\}$ note that the same analysis gives also for any $l \geq \rho n$, uniformly in $x \in S_{K_n}$, $\bar{\gamma}$ and $m_k \leq k!$,

$$
\begin{aligned}
(5.9) \quad & \mathbf{P}(N_{n,k}^x = m_k, k \in M_l) \\
& = (1 + o(1_n)) \prod_{k=l}^{n-1} \binom{m_{k+1} + m_k - 1}{m_k} p_k^{m_k} (1 - p_k)^{m_{k+1}}.
\end{aligned}
$$

Recall that $n_k(a) = 3ak^2 \log k$ and that we write $N \overset{k}{\sim} n_k$ if $|N - n_k| \leq k$ for $\rho n \leq k \leq n-1$ and $N = 0$ when $k = 0$. Relying upon the first moment estimates and Corollary 5.1, we next prove the second moment estimates (4.5) and (4.6). To this end, fix $x, y \in S_{K_n}$ with $2r_{l+1} + 2 > d(x, y) \geq 2r_l + 2$ for some $\rho'n \leq l \leq n-1$. Since $r_{l+2} - r_l \gg 2r_{l+1}$, it is easy to see that $D(y, r_l) \cap \partial D(x, r_k) = \varnothing$ for all $k \neq l+1$. Replacing hereafter $l$ by $l \wedge (n-3)$, it follows that for $k \neq l+1$, $k \neq l+2$, the events $\{N_{n,k}^x \overset{k}{\sim} n_k\}$ are measurable on the $\sigma$-algebra $\mathcal{G}_l^y$. With $J_l := \{0, \rho n, \ldots, l-1\}$ and $I_l := \{0, \rho n, \ldots, l, l+3, \ldots, n-1\}$, we note that

$$\{x, y \text{ are } n\text{-successful}\} \subset \{N_{n,k}^x \overset{k}{\sim} n_k, k \in I_l\} \cap \{N_{n,k}^y \overset{k}{\sim} n_k, k \in J_{l+1}\}.$$

Applying (5.1), we have that for some universal constant $C_3 < \infty$,

$$
\begin{aligned}
(5.10) \quad & \mathbf{P}(x \text{ and } y \text{ are } n\text{-successful}) \\
& \leq \sum_{m_l \overset{l}{\sim} n_l} \mathbb{E}[\mathbf{P}(N_{n,k}^y \overset{k}{\sim} n_k, k \in J_l | N_{n,l}^y = m_l, \mathcal{G}_l^y); N_{n,k}^x \overset{k}{\sim} n_k, k \in I_l] \\
& \leq C_3 \mathbf{P}(N_{n,k}^x \overset{k}{\sim} n_k, k \in I_l) \sum_{m_l \overset{l}{\sim} n_l} \mathbf{P}(N_{n,k}^y \overset{k}{\sim} n_k, k \in J_l | N_{n,l}^y = m_l).
\end{aligned}
$$

Using Corollary 5.1 once more, we have that

$$
\begin{aligned}
& (1 + o(1_n))\bar{q}_n \\
& = \mathbf{P}(y \text{ is } n\text{-successful}) \\
& = \sum_{m_l \overset{l}{\sim} n_l} \mathbb{E}[\mathbf{P}(N_{n,k}^y \overset{k}{\sim} n_k, k \in J_l | N_{n,l}^y = m_l, \mathcal{G}_l^y);
\end{aligned}
$$

(5.11)
$$N_{n,l}^y = m_l, N_{n,k}^y \overset{k}{\sim} n_k, k \in M_{l+1}]$$

$$\geq C_4 \sum_{m_l \overset{l}{\sim} n_l} \mathbf{P}(N_{n,l}^y = m_l, N_{n,k}^y \overset{k}{\sim} n_k, k \in M_{l+1})$$

$$\times \mathbf{P}(N_{n,k}^y \overset{k}{\sim} n_k, k \in J_l | N_{n,l}^y = m_l),$$

for some universal constant $C_4 > 0$. Hence by (5.9) and (5.7), for some universal constant $C_5 < \infty$,

(5.12) $$\sum_{m_l \overset{l}{\sim} n_l} \mathbf{P}(N_{n,k}^y \overset{k}{\sim} n_k, k \in J_l | N_{n,l}^y = m_l) \leq C_5^{n-l} l \left( \prod_{k=l}^{n-1} k^{3a} \sqrt{\log k} \right) \bar{q}_n.$$

Similarly, using Corollary 5.1,

$$\mathbf{P}(N_{n,k}^x \overset{k}{\sim} n_k, k \in I_l)$$

(5.13)
$$\leq \sum_{m_l \overset{l}{\sim} n_l} \mathbb{E}[\mathbf{P}(N_{n,k}^x \overset{k}{\sim} n_k, k \in J_l | N_{n,l}^x = m_l, \mathcal{G}_l^x); N_{n,k}^x \overset{k}{\sim} n_k, k \in M_{l+3}]$$

$$\leq C_6 \mathbf{P}(N_{n,k}^x \overset{k}{\sim} n_k, k \in M_{l+3}) \sum_{m_l \overset{l}{\sim} n_l} \mathbf{P}(N_{n,k}^x \overset{k}{\sim} n_k, k \in J_l | N_{n,l}^x = m_l).$$

Comparing (5.13) and (5.11), and applying once more (5.9) and (5.7), we get that

(5.14) $$\mathbf{P}(N_{n,k}^x \overset{k}{\sim} n_k, k \in I_l) \leq C_7 l \left( \prod_{k=l}^{l+2} k^{3a} \sqrt{\log k} \right) \bar{q}_n.$$

Putting (5.10), (5.12) and (5.14) together proves (4.5).

In case $d(x,y) \geq 2(r_n + 1)$, the event $\{x$ is $n$-successful$\}$ is $\mathcal{G}_n^y$ measurable, hence

$$\mathbf{P}(x \text{ and } y \text{ are } n\text{-successful})$$

$$= \mathbb{E}(\{\mathbf{P}(y \text{ is } n\text{-successful} | \mathcal{G}_n^y)\}; x \text{ is } n\text{-successful})$$

$$= \mathbb{E}(\{\mathbf{P}(N_{n,k}^y \overset{k}{\sim} n_k, k \in J_n | N_{n,n}^y = n_n, \mathcal{G}_n^y)\}; x \text{ is } n\text{-successful}),$$

and (4.6) follows from Corollary 5.1. $\square$

## 6. Large deviation bounds.

This section provides crucial large deviations estimates that are key to the proofs of Theorem 1.2 and of the upper bounds in Theorems 1.3 and 1.4. Roughly, we will be providing precise decay rates for the events that certain normalized excursion counts of balls concentric to a point $z$ (excursions between levels $r_{\beta n-1}$ and $r_{\beta n}$, before making $n_n$ excursions between levels $r_{n-1}$ and $r_n$ are atypical, together with forcing one or two points nearby not to be visited during these excursions.

More precisely, fix $0 < \beta < 1$ and $\tilde{n} \geq n$. Recall the definition $F_{h,\beta}(\gamma) = (1 - \gamma\beta)^2/(1 - \beta) + h\gamma^2\beta$ of (1.10). For any $h \geq 0$, the unique global minimum of $F_{h,\beta}(\gamma)$ is at $\gamma_h = \gamma_h(\beta) = 1/(h(1 - \beta) + \beta)$. For $0 < a < 2$, with $N_{n,k}^x = N_{n,k}^x(a)$ and $\mathcal{R}_n^x = \mathcal{R}_n^x(a)$ as in Section 4, we establish large deviations bounds away from $\gamma_h$ for the random variables $\widehat{N}_{n,\beta n}^x(a) := N_{n,\beta n}^x(a)/n_{\beta n}(a)$ together with the events $\{\mathcal{T}_{K_{\tilde{n}}}(x) > \mathcal{R}_n^z(a)\}$ and $\{\mathcal{T}_{K_{\tilde{n}}}(x') > \mathcal{R}_n^z(a)\}$ for $\tilde{n} \geq n$ and $x, x'$ not too far from $z \in \mathbb{Z}_{K_{\tilde{n}}}^2$, that is, $(z, x)$ and $(z, x, x')$ belonging to the sets

$$(6.1) \qquad \mathbf{Go}(\tilde{n}) = \{(z, x) : z \in \mathbb{Z}_{K_{\tilde{n}}}^2, x \in D(z, r_{\beta n - 2}) \cap \mathbb{Z}_{K_{\tilde{n}}}^2\},$$

$$(6.2) \qquad \begin{aligned} \mathbf{Go}^h(\tilde{n}) = \{(z, x, x') : (z, x) \in \mathbf{Go}(\tilde{n}), \\ (z, x') \in \mathbf{Go}(\tilde{n}), x' \notin D(x, r_{\beta h n/2 - 3})\}, \end{aligned}$$

where $h \in (0, 2)$. To express the bounds, define

$$I_h(\gamma) = \begin{cases} [0, \gamma^2], & \gamma < \gamma_h, \\ [0, \infty), & \gamma = \gamma_h, \\ [\gamma^2, \infty), & \gamma > \gamma_h. \end{cases}$$

LEMMA 6.1. *Fixing $0 < h < 2$ and $a, \gamma, \delta > 0$, for all $\tilde{n} \geq n \geq n_0$ we have the bounds*

$$(6.3) \qquad \max_{z \in \mathbb{Z}_{K_{\tilde{n}}}^2} \mathbf{P}(\widehat{N}_{n,\beta n}^z(a) \in I_0(\gamma)) \leq K_n^{-aF_{0,\beta}(\gamma) + \delta},$$

$$(6.4) \qquad \begin{aligned} \max_{(z,x) \in \mathbf{Go}(\tilde{n})} \mathbf{P}(\mathcal{T}_{K_{\tilde{n}}}(x) > \mathcal{R}_n^z(a), \widehat{N}_{n,\beta n}^z(a) \in I_1(\gamma)) \\ \leq K_n^{-aF_{1,\beta}(\gamma) + \delta}, \end{aligned}$$

$$(6.5) \qquad \begin{aligned} \max_{(z,x,x') \in \mathbf{Go}^h(\tilde{n})} \mathbf{P}(\mathcal{T}_{K_{\tilde{n}}}(x') > \mathcal{R}_n^z(a), \mathcal{T}_{K_{\tilde{n}}}(x) > \mathcal{R}_n^z(a), \widehat{N}_{n,\beta n}^z(a) \in I_h(\gamma)) \\ \leq K_n^{-aF_{h,\beta}(\gamma) + \delta}. \end{aligned}$$

As is often the case with large deviation statements, the key to the proof of Lemma 6.1 lies in the evaluation of certain moment generating functions. To state these, fix $z \in \mathbb{Z}_{K_{\tilde{n}}}^2$, and abbreviate $\partial_k$ for $\partial D(z, r_k)$. Consider a path of the simple random walk starting at a fixed $y \in \partial_{n-1}$. Let $Z$ denote the number of excursions of the path from $\partial_{\beta n - 1}$ to $\partial_{\beta n}$ until $T_{\partial_n}$ and $A(x) = \{T_{\partial_h}(x) < T_x\}$. See Figure 2.

Let $\lambda_h^* = 1/(1 - \beta) + h/\beta$ for $0 \leq h < 2$.

LEMMA 6.2. *Uniformly in $(z, x, x') \in \mathbf{Go}^h(\tilde{n})$, $\tilde{n} \geq n$ and $y \in \partial D(z, r_{n-1})$,*

$$(6.6) \qquad \mathbb{E}^y(e^{\lambda Z/n}) \leq 1 + \frac{1}{n}\left(\frac{\lambda}{1-(1-\beta)\lambda}\right) + \frac{c(0,\lambda)}{n\log n},$$

*for some $c(0,\lambda) < \infty$ and all $\lambda < \lambda_0^*$,*

$$(6.7) \qquad \mathbb{E}^y(e^{\lambda Z/n}\mathbf{1}_{A(x)}) \leq 1 + \frac{1}{n}\left(\frac{\beta\lambda-1}{\beta-(1-\beta)(\lambda\beta-1)}\right) + \frac{c(1,\lambda)}{n\log n},$$

*for some $c(1,\lambda) < \infty$ and all $\lambda < \lambda_1^*$, and*

$$(6.8) \quad \mathbb{E}^y(e^{\lambda Z/n}\mathbf{1}_{A(x)}\mathbf{1}_{A(x')}) \leq 1 + \frac{1}{n}\left(\frac{\beta\lambda-h}{\beta-(1-\beta)(\lambda\beta-h)}\right) + \frac{c(h,\lambda)}{n\log n},$$

*for some $c(h,\lambda) < \infty$ and all $\lambda < \lambda_h^*$.*

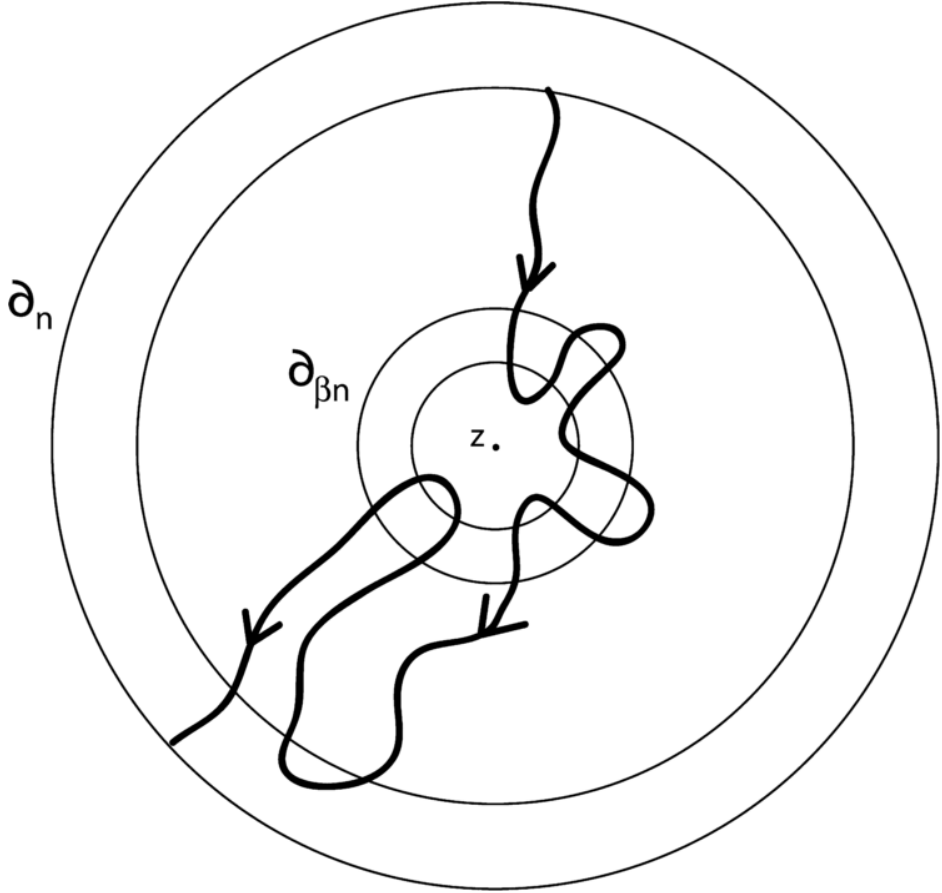

FIG. 2. $Z = 4$.

REMARK. The bound (6.8) is an improvement over (6.7), and will be used, only in the region $h > 1$ (in fact, $h$ near 2).

PROOF OF LEMMA 6.2. Recall that by (2.4), for some $c_1 < \infty$, all $n \geq n_0$ and any $z$:

$$(6.9) \quad q_- \leq \min_{v \in \partial_{\beta n}} \mathbf{P}^v(T_{\partial_n} < T_{\partial_{\beta n-1}}) \leq \max_{v \in \partial_{\beta n}} \mathbf{P}^v(T_{\partial_n} < T_{\partial_{\beta n-1}}) \leq q_+,$$

$$(6.10) \quad q_- \leq \min_{v \in \partial_{n-1}} \mathbf{P}^v(T_{\partial_{\beta n-1}} < T_{\partial_n}) \leq \max_{v \in \partial_{n-1}} \mathbf{P}^v(T_{\partial_{\beta n-1}} < T_{\partial_n}) \leq q_+,$$

where $q_\pm = (1 - \beta)^{-1} n^{-1}(1 \pm c_1/\log n)$. By (6.10), for any $y \in \partial_{n-1}$,

$$(6.11) \quad \mathbf{P}^y(Z = 0) = \mathbf{P}^y(T_{\partial_{\beta n-1}} > T_{\partial_n}) \leq 1 - q_-,$$

and for $j = 1, 2, \ldots$ we have $Z = j$ if we first visit $\partial_{\beta n-1}$ prior to $\partial_n$, then have exactly $j - 1$ cycles consisting of visits to $\partial_{\beta n}$ and back to $\partial_{\beta n-1}$, prior to the first visit to $\partial_n$. Hence, by (6.9), (6.10) and the strong Markov property, for any $y \in \partial_{n-1}$ we have that $\mathbf{P}^y(Z = j) \leq (1 - q_-)^{j-1} q_+^2$. The bound (6.6) then follows from the $h = 0$ case of the inequality

$$(6.12) \quad \begin{aligned} &(1 - q_-) + \sum_{j=1}^{\infty} e^{\lambda j/n}(1 - p_h)^j (1 - q_-)^{j-1} q_+^2 \\ &\leq 1 + \frac{1}{n}\left(\frac{\beta\lambda - h}{\beta - (1 - \beta)(\lambda\beta - h)}\right) + \frac{c(h, \lambda)}{n \log n}, \end{aligned}$$

where in general

$$p_h = \frac{h}{\beta n}(1 - c'/\log n) \quad \text{and} \quad \lambda < \lambda_h^*.$$

To see (6.12) let $v = 1/(1 - \beta) + h/\beta - \lambda$. Then, for some finite $C$ and $n_0$ (both depending on $c'$, $c_1$, $\lambda$, $h$ and $\beta$), we have that $q_+^2(1 - p_h)e^{\lambda/n} \leq (n(1 - \beta))^{-2}(1 + C/\log n)$ and $1 - e^{\lambda/n}(1 - p_h)(1 - q_-) \geq n^{-1}v(1 - C/\log n)$, for all $n \geq n_0$. Consequently, for some $c = c(h, \lambda) < \infty$ and all $n \geq n_0$,

$$\begin{aligned} &(1 - q_-) + \sum_{j=1}^{\infty} e^{\lambda j/n}(1 - p_h)^j(1 - q_-)^{j-1} q_+^2 \\ &= (1 - q_-) + \frac{q_+^2(1 - p_h)e^{\lambda/n}}{1 - e^{\lambda/n}(1 - p_h)(1 - q_-)} \\ &\leq 1 + \frac{1}{n}\left(\frac{1}{(1 - \beta)^2 v} - \frac{1}{1 - \beta}\right) + \frac{c}{n \log n} \\ &= 1 + \frac{1}{n}\left(\frac{\beta\lambda - h}{\beta + h(1 - \beta) - \beta(1 - \beta)\lambda}\right) + \frac{c}{n \log n}, \end{aligned}$$

which gives (6.12).

We next turn to (6.7). Enlarging $c_1$ as needed, by (2.1) we have that for all $\tilde{n} \geq n \geq n_0$ and $(z, x) \in \mathbf{Go}(\tilde{n})$,

$$
\begin{aligned}
\min_{v \in \partial_{\beta_{n-1}}} \mathbf{P}^v(T_x < T_{\partial_{\beta_n}}) &\geq \min_{v \in D(x, 2r_{\beta_{n-1}})} \mathbf{P}^v(T_x < T_{\partial D(x, 0.5 r_{\beta_n})}) \\
&\geq \frac{1}{\beta n}\left(1 - \frac{c_1}{\log n}\right) =: p.
\end{aligned}
$$

(6.13)

We have $Z \mathbf{1}_{A(x)} = j \geq 1$ if we first visit $\partial_{\beta_{n-1}}$ prior to $\partial_n$, then have $j - 1$ cycles consisting of visits to $\partial_{\beta_n}$ and back to $\partial_{\beta_{n-1}}$ without hitting $x$ or $\partial_n$, and finally, a visit to $\partial_n$ without hitting $x$. Hence, by (6.9), (6.10), (6.13) and the strong Markov property, for any $z$, $y$ and $x$ as above,

$$
\mathbf{P}^y(Z = j, A(x)) \leq (1 - p)^j (1 - q_-)^{j-1} q_+^2.
$$

(6.14)

Note that $A(x)$ occurs when $Z = 0$, so that (6.11), (6.14) and the $h = 1$ case of (6.12) give (6.7).

We finally turn to (6.8). By the strong Markov property at $\min(T_x, T_{x'})$, for $v \in \partial_{\beta_{n-1}}$ and $x, x' \in D(z, r_{\beta_{n-2}})$,

$$
\begin{aligned}
&\mathbf{P}^v(\max(T_x, T_{x'}) < T_{\partial_{\beta_n}}) \\
&\leq \mathbf{P}^v(T_{x'} < T_{\partial_{\beta_n}}) \mathbf{P}^{x'}(T_x < T_{\partial_{\beta_n}}) \\
&\quad + \mathbf{P}^v(T_x < T_{\partial_{\beta_n}}) \mathbf{P}^x(T_{x'} < T_{\partial_{\beta_n}}).
\end{aligned}
$$

(6.15)

Enlarging $c_1$ as needed, since $\log r_{\beta_{hn/2-3}} / \log r_{\beta_n} = h/2 + O(1/\log n)$, similarly to the derivation of (6.13) we have by (2.1) that for all $n \geq n_0$ and $(z, x, x') \in \mathbf{Go}^h(\tilde{n})$,

$$
\begin{aligned}
&\mathbf{P}^x(T_{x'} < T_{\partial_{\beta_n}}) \max_{v \in \partial_{\beta_{n-1}}} \mathbf{P}^v(T_x < T_{\partial_{\beta_n}}) \\
&\leq \mathbf{P}^x(T_{x'} < T_{\partial D(x', 2r_{\beta_n})}) \max_{d(v, x) \geq 0.5 r_{\beta_{n-1}}} \mathbf{P}^v(T_x < T_{\partial D(x, 2r_{\beta_n})}) \\
&\leq \frac{1}{\beta n}\left(1 - \frac{h}{2} + \frac{c_1}{\log n}\right).
\end{aligned}
$$

The same bound applies to the other term on the right-hand side of (6.15). When combined with (6.13) which applies for both $x$ and $x'$, these bounds yield (by inclusion–exclusion) that for all $n \geq n_0$, uniformly in $(z, x, x') \in \mathbf{Go}^h(\tilde{n})$,

$$
\begin{aligned}
&\max_{v \in \partial_{\beta_{n-1}}} \mathbf{P}^v(T_x > T_{\partial_{\beta_n}}, T_{x'} > T_{\partial_{\beta_n}}) \\
&\leq 1 - 2p + \frac{2}{\beta n}\left(1 - \frac{h}{2} + \frac{c_1}{\log n}\right) =: 1 - \widehat{p}_h
\end{aligned}
$$

(6.16)

with $\widehat{p}_h = \frac{h}{\beta n}(1 - c'/\log n)$ of the same form as $p_h$. Note that $Z\mathbf{1}_{A(x)}\mathbf{1}_{A(x')} = j \geq 1$ if the walk visits $\partial_{\beta n - 1}$ prior to $\partial_n$, then has $j-1$ cycles consisting of visits to $\partial_{\beta n}$ and back to $\partial_{\beta n - 1}$, without hitting $x$, $x'$ or $\partial_n$, and finally, visits to $\partial_n$ without hitting $x$ or $x'$. Hence, by (6.9), (6.10), (6.16) and the strong Markov property, for any $z$, $y$, $x$ and $x'$ as above,

$$\mathbf{P}^y(Z = j, A(x), A(x')) \leq (1 - \widehat{p}_h)^j (1 - q_-)^{j-1} q_+^2,$$

and (6.8) now follows as in the derivation of (6.7). This completes the proof of Lemma 6.2. $\quad\square$

PROOF OF LEMMA 6.1. A straightforward calculation shows that for any $h \geq 0$ and $\gamma > 0$,

$$F_{h,\beta}(\gamma) = \lambda_{h,\gamma} \gamma^2 \beta^2 - \frac{\beta \lambda_{h,\gamma} - h}{\beta - (1 - \beta)(\lambda_{h,\gamma}\beta - h)}$$

(6.17)
$$\text{where } \lambda_{h,\gamma} := \frac{\beta + h(1 - \beta) - 1/\gamma}{\beta(1 - \beta)} < \lambda_h^*,$$

and $\lambda_{h,\gamma} \leq 0$ if and only if $\gamma \leq \gamma_h$.

Let $\widehat{Z}_0$ denote the number of excursions from $\partial_{\beta n - 1}$ to $\partial_{\beta n}$ before $X_t$ first hits $\partial_{n-1}$ and let $A_0(x)$ denote the event that $x$ is not visited during this time interval. For any $j \geq 1$ let $\widehat{Z}_j$ denote the number of excursions from $\partial_{\beta n - 1}$ to $\partial_{\beta n}$ during the $j$th excursion of $X_t$ from $\partial_{n-1}$ to $\partial_n$ and let $A_j(x)$ denote the event that $x$ is not visited during this excursion. With this notation,

$$N_{n,\beta n}^z(a) = \sum_{j=0}^{3an^2 \log n} \widehat{Z}_j,$$

and the event $\{\mathcal{T}_{K_n}(x) > \mathcal{R}_n^z(a)\}$ is the intersection of the events $A_j(x)$ for $j = 0, \ldots, 3an^2 \log n$. Consequently, using Chebyshev's inequality and the strong Markov property (at the start of the $3an^2 \log n$ excursions from $\partial_{n-1}$ to $\partial_n$), for any $\delta > 0 \geq \lambda$ and all $\tilde{n} \geq n \geq n_0$, uniformly in $z$,

$$\mathbf{P}(\widehat{N}_{n,\beta n}^z(a) \leq \gamma^2) \leq e^{-\lambda \gamma^2 n_{\beta n}/n} \mathbb{E}(e^{(\lambda/n)\sum_{j=1}^{n_n} \widehat{Z}_j})$$

(6.18)
$$\leq K_n^{-a\lambda\gamma^2\beta^2 + \delta} \left( \max_{y \in \partial_{n-1}} \mathbb{E}^y(e^{\lambda Z/n}) \right)^{3an^2 \log n}.$$

Per $\gamma \leq \gamma_0$ consider (6.18) for $\lambda = \lambda_{0,\gamma} \leq 0$, applying (6.17) and (6.6) to obtain (6.3) in case $\gamma < \gamma_0$. Turning to deal with $\gamma \geq \gamma_0$, note that $\mathbf{P}^y(Z = j) \leq (1 - q_-)^{j-1} q_+$ for all $j \geq 1$, even if $y \in \partial_{\beta n - 1}$. Thus, for any $\lambda < \lambda_0^*$, similar to the derivation of (6.6) we get that for some $c_5 = c_5(\lambda) < \infty$ and all $\tilde{n} \geq n \geq n_0$, uniformly in $z$,

$$(6.19) \qquad \mathbb{E}(e^{\lambda \widehat{Z}_0/n}) \leq \max_{y \in \partial_{\beta n - 1}} \mathbb{E}^y(e^{\lambda Z/n}) \leq c_5.$$

Analogous to (6.18) we also have that for any $\delta > 0$, $\lambda \geq 0$, $\tilde{n} \geq n \geq n_0$ and $z$,

$$(6.20) \qquad \begin{aligned} \mathbf{P}(\widehat{N}_{n,\beta n}^z(a) \geq \gamma^2) &\leq e^{-\lambda \gamma^2 n_{\beta n}/n} \mathbb{E}(e^{(\lambda/n)\sum_{j=0}^{n_n} \widehat{Z}_j}) \\ &\leq c_5 K_n^{-a\lambda \gamma^2 \beta^2 + \delta} \left( \max_{y \in \partial_{n-1}} \mathbb{E}^y(e^{\lambda Z/n}) \right)^{n_n}. \end{aligned}$$

Considering (6.20) for $\lambda = \lambda_{0,\gamma} \geq 0$ (as $\gamma \geq \gamma_0$), and applying (6.6) and (6.17), we complete the proof of (6.3).

Similarly, we have that for any $\delta > 0 \geq \lambda$, $\tilde{n} \geq n \geq n_0$, $(z,x) \in \mathbf{Go}(\tilde{n})$,

$$(6.21) \qquad \begin{aligned} &\mathbf{P}(\mathcal{T}_{K_{\tilde{n}}}(x) > \mathcal{R}_n^z(a), \widehat{N}_{n,\beta n}^z(a) \leq \gamma^2) \\ &\quad \leq e^{-\lambda \gamma^2 n_{\beta n}/n} \mathbb{E}\left( \prod_{j=1}^{n_n} e^{(\lambda/n)\widehat{Z}_j} \mathbf{1}_{A_j(x)} \right) \\ &\quad \leq K_n^{-a\lambda \gamma^2 \beta^2 + \delta} \left( \max_{y \in \partial_{n-1}} \mathbb{E}^y(e^{\lambda Z/n} \mathbf{1}_{A(x)}) \right)^{3an^2 \log n}. \end{aligned}$$

Given $\gamma \leq \gamma_1$, consider (6.21) for $\lambda = \lambda_{1,\gamma} \leq 0$, and apply (6.17) and (6.7) to get (6.4) for $\gamma < \gamma_1$. Further, the same argument leading to (6.19) shows also that $\max_{y \in \partial_{\beta n - 1}} \mathbb{E}^y(e^{\lambda Z/n} \mathbf{1}_{A(x)}) \leq c_5$ for all $\lambda < \lambda_1^*$. Consequently, for $\delta > 0$, $\lambda \geq 0$, $\tilde{n} \geq n \geq n_0$ and $(z,x) \in \mathbf{Go}(\tilde{n})$,

$$\begin{aligned} &\mathbf{P}(\mathcal{T}_{K_{\tilde{n}}}(x) > \mathcal{R}_n^z(a), \widehat{N}_{n,\beta n}^z(a) \geq \gamma^2) \\ &\quad \leq c_5 K_n^{-a\lambda \gamma^2 \beta^2 + \delta} \left( \max_{y \in \partial_{n-1}} \mathbb{E}^y(e^{\lambda Z/n} \mathbf{1}_{A(x)}) \right)^{n_n}, \end{aligned}$$

and since $\lambda_{1,\gamma} \geq 0$ for $\gamma \geq \gamma_1$, we complete the proof of (6.4) by using again (6.17) and (6.7).

Using (6.8) and $\lambda = \lambda_{h,\gamma}$, the proof of (6.5) proceeds along the same lines, thus completing the proof of Lemma 6.1. $\square$

## 7. Late points in a small neighborhood. 
We devote this section to the proof of Theorem 1.2, as the basic large deviations bounds needed are already in place.

PROOF OF THEOREM 1.2. We actually show that for $0 < \alpha < \beta^2 < 1$, some $b < \infty$, any $\xi, \delta, \eta > 0$, and all $n \geq n_0$, $\bar{\gamma} \in \mathcal{I}$ and $x = x_n \in \mathbb{Z}_{K_n}^2$,

$$(7.1) \qquad \mathbf{P}(|\mathcal{L}_{K_n}(\alpha) \cap D(x, r_{\beta n - b})| \geq K_n^{2\beta - (2\alpha - \xi)/\beta + 4\delta}) \leq 2\eta,$$

(7.2)          $\mathbf{P}(|\mathcal{L}_{K_n}(\alpha) \cap D(x, r_{\beta n + b})| \geq K_n^{2\beta - (2\alpha + \xi)/\beta - \delta}) \geq 1 - 2\eta.$

Since $\log r_{\beta n \pm b}/\log K_n \to \beta$ and the set of $K_n$ values cover all large integers, the theorem follows by considering $\eta \downarrow 0$ and adjusting the values of $\beta$, $\delta > 0$ and $\xi > 0$.

Starting with the upper bound (7.1), recall the notation $\mathcal{R}_k^x(a)$ for the time until completion of the first $n_k(a) = 3ak^2 \log k$ excursions from $\partial D(x, r_{k-1})$ to $\partial D(x, r_k)$, $k = 3, \ldots, n$, then $N_{k,0}^x(a)$ for the number of visits to $x$ until time $\mathcal{R}_k^x(a)$, and $N_{k,l}^x(a)$, $2 \leq l \leq k-1$, for the number of excursions from $\partial D(x, r_{l-1})$ to $\partial D(x, r_l)$ until time $\mathcal{R}_k^x(a)$. Let $t_n^* = \frac{4}{\pi}(K_n \log K_n)^2$ and

(7.3)          $\widehat{\mathcal{L}}_{K_n}(\tilde{a}) := \left\{ y \in \mathbb{Z}_{K_n}^2 : \mathcal{T}_{K_n}(y) > \max_{z \in \mathbb{Z}_{K_n}^2} \mathcal{R}_n^z(\tilde{a}) \right\},$

taking hereafter $\xi \in (0, 2\alpha)$ and $\tilde{a} = 2\alpha - \xi > 0$ (in the remainder of the paper we always have $\tilde{a} < 2\alpha < a$). Applying (3.19) with $R = r_n$, $r = r_{n-1}$ and $N = 3\tilde{a}n^2 \log n$, we see that for some $c = c(\alpha, \xi) > 0$ and all $n$,

(7.4)          $\max_{z \in \mathbb{Z}_{K_n}^2} \mathbf{P}(\mathcal{R}_n^z(\tilde{a}) \geq \alpha t_n^*) \leq c^{-1} e^{-cn^2 \log n},$

resulting with

(7.5)          $\lim_{n \to \infty} \mathbf{P}(\mathcal{L}_{K_n}(\alpha) \subseteq \widehat{\mathcal{L}}_{K_n}(\tilde{a})) = 1.$

Hence, to establish (7.1) it suffices to show that

(7.6)          $\mathbf{P}(|\widehat{\mathcal{L}}_{K_n}(\tilde{a}) \cap D(x, r_{\beta n - 2})| \geq K_n^{2\beta - \tilde{a}/\beta + 4\delta}) \leq \eta.$

Since $F_{0,\beta}(\gamma) > 0$ for $\gamma < \gamma_0 = 1/\beta$, it follows from (6.3) that for any $\delta' > 0$,

(7.7)          $\lim_{n \to \infty} \max_{x \in \mathbb{Z}_{K_n}^2} \mathbf{P}(N_{n,\beta n}^x(\tilde{a}) < (1 - \delta')n_n(\tilde{a})) = 0.$

Recall that $F_{1,\beta}(1/\beta) = 1/\beta$ and $r_{\beta n} \leq K_n^\beta$ for all $n$. Moreover, $(1 - \delta')n_n \geq \gamma^2 n_{\beta n}$ for $\gamma = (1 - \delta')/\beta$ and all $n$. Hence, if $\gamma \geq \gamma_1$, then by (6.4) we have that

$\mathbf{P}(|\widehat{\mathcal{L}}_{K_n}(\tilde{a}) \cap D(x, r_{\beta n - 2})| \geq K_n^{2\beta - \tilde{a}/\beta + 4\delta}, N_{n,\beta n}^x(\tilde{a}) \geq (1 - \delta')n_n(\tilde{a}))$

(7.8)          $\leq K_n^{-(2\beta - \tilde{a}/\beta) - 4\delta} r_{\beta n}^2 \max_{y \in D(x, r_{\beta n - 2})} \mathbf{P}(\mathcal{T}_{K_n}(y) > \mathcal{R}_n^x(\tilde{a}), \widehat{N}_{n,\beta n}^x(\tilde{a}) \geq \gamma^2)$

$\leq K_n^{\tilde{a}(F_{1,\beta}(1/\beta) - F_{1,\beta}(\gamma)) - 3\delta}.$

With $\beta < 1$, for $\delta' > 0$ small enough we have both $\gamma \geq \gamma_1 = 1$ and $F_{1,\beta}(\frac{1}{\beta}) - F_{1,\beta}(\gamma) \leq \delta$. Thus, considering (7.7) and (7.8) for such $\delta'$ completes the proof of (7.6), hence also of (7.1).

Turning to prove the lower bound (7.2), fixing $0 < \xi < 2(\beta^2 - \alpha)$ so $a' = (2\alpha + \xi)/\beta^2 < 2$ and $0 < \rho < (2 - a')/2$ we say that a point $y \in \mathbb{Z}^2_{K_n}$ is $\beta n$-*successful* if

$$N^y_{\beta n, 0}(a') = 0, \qquad N^y_{\beta n, k}(a') \overset{k}{\sim} n_k(a') \qquad \forall k = \rho \beta n, \dots, \beta n - 1.$$

In particular, if $y$ is $\beta n$-successful, then $\mathcal{T}_{K_n}(y) > \mathcal{R}^y_{\beta n}(a')$. Let $\mathcal{L}^\sharp_{K_n}(a', \beta n)$ be the set of points in $\mathbb{Z}^2_{K_n}$ which are $\beta n$-successful. A rerun of the proof of (4.3), this time with $\beta n$ replacing $n$, shows that for some $b \geq 10$, any $\delta > 0$, $\eta > 0$, all $n \geq n_0$, $\bar\gamma \in \mathcal{I}$ and $x \in \mathbb{Z}^2_{K_n}$,

$$(7.9) \qquad \mathbf{P}(|\mathcal{L}^\sharp_{K_n}(a', \beta n) \cap D(x, r_{\beta n + b})| \geq K_n^{\beta(2-a')-\delta}) \geq 1 - \eta.$$

Consequently, (7.2) follows once we show that uniformly in $x$,

$$(7.10) \qquad \mathbf{P}\left(\min_{y \in D(x, r_{\beta n + b})} \mathcal{R}^y_{\beta n}(a') \leq \alpha t^*_n\right) \to 0.$$

To see this, let $\mathcal{Y}_n$ be a minimal set of points in $D(x, r_{\beta n + b})$ such that

$$D(x, r_{\beta n + b}) \subseteq \bigcup_{y \in \mathcal{Y}_n} D(y, r_{\beta n - 2}).$$

Let $\widehat{\mathcal{R}}^y_{\beta n}(a')$ denote the time until completion of the first $n_{\beta n}(a')$ excursions from $\partial D(y, r_{\beta n - 1} + r_{\beta n - 2})$ to $\partial D(y, r_{\beta n} - r_{\beta n - 2})$. For any $z \in D(y, r_{\beta n - 2})$ we have that

$$D(z, r_{\beta n - 1}) \subseteq D(y, r_{\beta n - 1} + r_{\beta n - 2}) \subseteq D(y, r_{\beta n} - r_{\beta n - 2}) \subseteq D(z, r_{\beta n}),$$

implying that each excursion from $\partial D(z, r_{\beta n - 1})$ to $\partial D(z, r_{\beta n})$ requires at least one excursion from $\partial D(y, r_{\beta n - 1} + r_{\beta n - 2})$ to $\partial D(y, r_{\beta n} - r_{\beta n - 2})$. See Figure 3.

Thus, $\mathcal{R}^z_{\beta n}(a') \geq \widehat{\mathcal{R}}^y_{\beta n}(a')$ and consequently,

$$(7.11) \qquad \mathbf{P}\left(\min_{z \in D(y, r_{\beta n - 2})} \mathcal{R}^z_{\beta n}(a') \leq \alpha t^*_n\right) \leq \mathbf{P}(\widehat{\mathcal{R}}^y_{\beta n}(a') \leq \alpha t^*_n).$$

Applying (3.18) with $R = r_{\beta n} - r_{\beta n - 2}$, $r = r_{\beta n - 1} + r_{\beta n - 2}$ and $N = n_{\beta n}(a') = 3(2\alpha + \xi)n^2 \log(\beta n)$, the right-hand side of (7.11) is bounded by

$$C^{-1} \exp\left\{-C\left(\frac{\log(R/r)}{\log(K_n/r)}\right)n^2 \log n\right\} \leq c^{-1} \exp\{-cn \log n\},$$

for some $C, c > 0$ that depend only on $\alpha, \xi > 0$, yielding (7.10) (recall that $|\mathcal{Y}_n| \leq C n^{10b}$). $\quad\square$

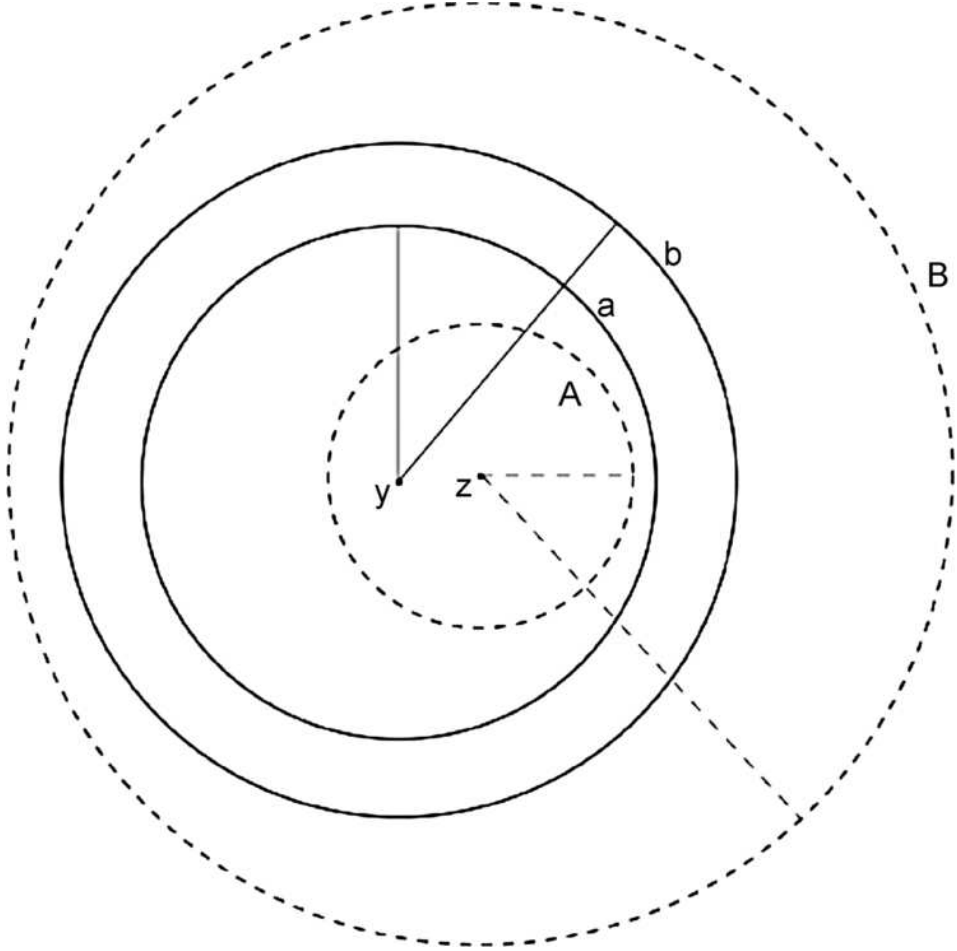

Fig. 3. $A = \partial D(z, r_{\beta_{n-1}}), a = \partial D(y, r_{\beta_{n-1}} + r_{\beta_{n-2}}), b = \partial D(y, r_{\beta_n} - r_{\beta_{n-2}}), B = \partial D(z, r_{\beta_n})$.

**8. Clusters of late points.** Fixing $0 < \alpha, \beta < 1$, this section is devoted to the proof of Theorem 1.3. As usual, it suffices to establish (1.5) and (1.6) for the subsequence $K_n = n^{\bar{\gamma}}(n!)^3$, provided all our estimates are uniform in $\bar{\gamma} \in \mathcal{I}$. To this end, set

$$(8.1) \qquad W^x(\beta_2, \beta_1) = |\{y \in \mathcal{L}_{K_n}(\alpha) : r_{\beta_2 n - 3} < d(x, y) \le r_{\beta_1 n - 3}\}|,$$

with $W^x = W^x(0, \beta)$. We actually prove that:

LEMMA 8.1. *For each $\delta > 0$ there exists $\varepsilon \in (0, \delta/2)$ such that*

$$(8.2) \qquad p_n := K_n^{2\alpha + \varepsilon} \max_{x \in \mathbb{Z}_{K_n}^2} \mathbf{P}(x \in \mathcal{L}_{K_n}(\alpha), W^x \le K_n^{2\beta(1-\alpha) - 5\delta}) \underset{n \to \infty}{\longrightarrow} 0.$$

Lemma 8.2. *For each $\delta > 0$ there exists $\varepsilon \in (0, \delta/2)$ such that*

$$(8.3) \quad \bar{p}_n := K_n^{2\alpha+\varepsilon} \max_{x \in \mathbb{Z}_{K_n}^2} \mathbf{P}(x \in \mathcal{L}_{K_n}(\alpha), W^x \geq K_n^{2\beta(1-\alpha)+5\delta}) \xrightarrow[n\to\infty]{} 0.$$

By (1.2), we have $\mathbf{P}(|\mathcal{L}_{K_n}(\alpha)| \geq K_n^{2(1-\alpha)-\varepsilon/2}) \to 1$ for $n \to \infty$, and with $\log r_{\beta n-3}/\log K_n \to \beta$, the bounds (8.2), (8.3) imply that (1.6) holds (adjusting $\beta$ as needed). These bounds also imply that (1.5) is a consequence of the uniform lower bound $\mathbf{P}(x \in \mathcal{L}_{K_n}(\alpha)) \geq K_n^{-2\alpha-\varepsilon/2}$, holding for any $n$ large enough and all $x \in \mathbb{Z}_{K_n}^2$, $x \neq 0$. Applying Lemma 4.1 we get the latter bound as soon as

$$(8.4) \quad \min_{x \in \mathbb{Z}_{K_n}^2 \setminus \{0\}} \mathbf{P}(\mathcal{T}_{K_n}(x) > \mathcal{R}_n^x(a)) \geq K_n^{-2\alpha-\varepsilon/3}$$

holds for $a = 2\alpha + \varepsilon/7$ and all $n$ sufficiently large. Since $\mathcal{T}_{K_n}(x) > \mathcal{R}_n^x(a)$ whenever $x$ is $n$-successful, by (4.4) and translation invariance of the SRW we have that

$$(8.5) \quad \min_{x \in \mathbb{Z}_{K_n}^2} \min_{y \notin D(x, r_n)} \mathbf{P}^y(\mathcal{T}_{K_n}(x) > \mathcal{R}_n^x(a)) \geq K_n^{-2\alpha-\varepsilon/6}.$$

For any finite $r > 0$ there exists $c = c(r) > 0$ such that $\mathbf{P}(T_x > T_{\partial D(x,r)}) \geq c$ for all $n$ sufficiently large and all $x \neq 0$. Consequently, by (2.1) we have that $\mathbf{P}(T_x > T_{\partial D(x,r_n)}) \geq c'/\log r_n \geq K_n^{-\varepsilon/6}$ for some $c' > 0$, all $n$ sufficiently large and all $x \neq 0$. Combining this with (8.5) and the strong Markov property at $T_{\partial D(x,r_n)}$ results with (8.4), thus completing the proof of Theorem 1.3. $\square$

Proof of Lemma 8.1. Let $\widehat{\mathcal{Z}}_{n,\beta}^x := \{z \in \widehat{\mathcal{Z}}_{n,\beta} : z \neq 0, d(x,z) < 0.5 r_{\beta n-3}\}$, where $\widehat{\mathcal{Z}}_{n,\beta'}$ denotes for each $0 < \beta' < 1$ a subgrid of $\mathbb{Z}_{K_n}^2$ of spacing $4r_{\beta'n-4}$ such that $0 \in \widehat{\mathcal{Z}}_{n,\beta'}$. Fixing $\xi \in (0, 2\alpha)$ and $\eta \in (0,1)$ to be chosen later, let $a = 2\alpha + \xi$, $a' = (1+2\eta)^3 a$ and

$$\widetilde{W}^z = |\{y \in D(z, r_{\beta n-6}) : \mathcal{T}_{K_n}(y) > \mathcal{R}_{\beta n-4}^z(a')\}|.$$

Fixing $x \in \mathbb{Z}_{K_n}^2$, let $\mathbf{C} := \{W^x \leq K_n^{2\beta(1-\alpha)-5\delta}\}$ and for any $z \in \widehat{\mathcal{Z}}_{n,\beta}^x$ define the events $\mathbf{A}_z := \{\mathcal{R}_n^z(a) \geq a t_n^*\}$, $\mathbf{B}_z := \{N_{n,\beta n-4}^z(a) \leq n_{\beta n-4}(a')\}$ and $\mathbf{C}_z = \{\widetilde{W}^z \leq K_n^{2\beta(1-\alpha)-5\delta}\}$. Observe that $\mathbf{A}_z \cap \mathbf{B}_z$ implies that $W^x \geq \widetilde{W}^z$ and hence $\mathbf{A}_z \cap \mathbf{B}_z \cap \mathbf{C} \subseteq \mathbf{C}_z$ for any $z \in \widehat{\mathcal{Z}}_{n,\beta}^x$. Further, setting $\tilde{a} = 2\alpha - \xi$ and considering the events $\mathbf{F}_z := \{\mathcal{T}_{K_n}(x) > \mathcal{R}_n^z(\tilde{a})\}$ and $\mathbf{H}_z := \{\mathcal{R}_n^z(\tilde{a}) \geq a t_n^*\}$, we have that if $x \in \mathcal{L}_{K_n}(\alpha)$, then $\mathbf{H}_z \cup \mathbf{F}_z$ holds for each $z \in \widehat{\mathcal{Z}}_{n,\beta}^x$. Note that by the preceding $\mathbf{A}_z \cap \mathbf{F}_z \cap \mathbf{C} \subseteq (\mathbf{F}_z \cap \mathbf{B}_z^c) \cup \mathbf{C}_z$ for each $z \in \widehat{\mathcal{Z}}_{n,\beta}^x$ and hence

$$\{x \in \mathcal{L}_{K_n}(\alpha), \mathbf{C}\} \subseteq \bigcup_z \mathbf{H}_z \bigcup_z \mathbf{A}_z^c \cup \left(\bigcap_z (\mathbf{A}_z \cap \mathbf{F}_z \cap \mathbf{C})\right)$$

$$\subseteq \bigcup_z \mathbf{H}_z \bigcup_z \mathbf{A}_z^c \cup \left( \bigcap_z \mathbf{C}_z \right) \bigcup_z (\mathbf{F}_z \cap \mathbf{B}_z^c).$$

With $|\widehat{\mathcal{Z}}_{n,\beta}^x| \leq K_n^\varepsilon$ for all $\varepsilon > 0$ and $n$ sufficiently large, we thus have that

$$p_n \leq K_n^{2\alpha+\varepsilon} \mathbf{P}\left( \max_{z \in \mathbb{Z}_{K_n}^2} \mathcal{R}_n^z(\tilde{a}) \geq \alpha t_n^* \right) + K_n^{2\alpha+\varepsilon} \mathbf{P}\left( \min_{z \in \mathbb{Z}_{K_n}^2} \mathcal{R}_n^z(a) \leq \alpha t_n^* \right)$$

$$+ K_n^{2\alpha+\varepsilon} \max_{x \in \mathbb{Z}_{K_n}^2} \mathbf{P}\left( \max_{z \in \widehat{\mathcal{Z}}_{n,\beta}^x} \widetilde{W}^z \leq K_n^{2\beta(1-\alpha)-5\delta} \right)$$

$$+ K_n^{2\alpha+2\varepsilon} \max_{x,z \in \widehat{\mathcal{Z}}_{n,\beta}^x} \mathbf{P}(\mathcal{T}_{K_n}(x) > \mathcal{R}_n^z(\tilde{a}), N_{n,\beta n-4}^z(a) > n_{\beta n-4}(a'))$$

$$:= p_{n,0} + p_{n,1} + p_{n,2} + p_{n,3}.$$

By (7.4) we have that $p_{n,0} \to 0$ as $n \to \infty$. With $a > 2\alpha$, by (3.18), similar to the derivation of (7.4) we get also that $p_{n,1} \to 0$ as $n \to \infty$.

Turning to deal with the term $p_{n,2}$, consider the $\sigma$-algebra $\mathcal{G} = \bigcap_{z \in \widehat{\mathcal{Z}}_{n,\beta}^x} \mathcal{G}^z$, for $\mathcal{G}^z$ corresponding to $R' = r_{\beta n-4}$ and $R = r_{\beta n-5}$ in Lemma 2.4. Since $D(z', r_{\beta n-4}) \subseteq D(z, r_{n-1}) \setminus D(z, r_{\beta n-4})$ for any $z, z' \in \widehat{\mathcal{Z}}_{n,\beta}^x$, it follows that conditional upon $\mathcal{G}$, the random variables $\{\widetilde{W}^z\}_{z \in \widehat{\mathcal{Z}}_{n,\beta}^x}$ are independent with $\widetilde{W}^z$ measurable on the $\sigma$-algebra $\mathcal{H}^z(n_{\beta n-4}(a'))$ corresponding to $r = r_{\beta n-6}$ in Lemma 2.4. With $|\widehat{\mathcal{Z}}_{n,\beta}^x| \geq n^2$ for all $n$ sufficiently large, it follows from the latter lemma that

$$p_{n,2} = K_n^{2\alpha+\varepsilon} \max_{x \in \mathbb{Z}_{K_n}^2} \mathbb{E}\left( \prod_{z \in \widehat{\mathcal{Z}}_{n,\beta}^x} \mathbf{P}(\widetilde{W}^z \leq K_n^{2\beta(1-\alpha)-5\delta}|\mathcal{G}) \right)$$

$$\leq K_n^{2\alpha+\varepsilon} \left( (1+o(1_n)) \max_{z \in \mathbb{Z}_{K_n}^2} \mathbf{P}(\widetilde{W}^z \leq K_n^{2\beta(1-\alpha)-5\delta}) \right)^{n^2} \xrightarrow[n \to \infty]{} 0,$$

provided that for some universal constant $c > 0$

$$(8.6) \qquad \min_{z \in \mathbb{Z}_{K_n}^2} \mathbf{P}(\widetilde{W}^z \geq K_n^{2\beta(1-\alpha)-5\delta}) \geq c.$$

Applying (3.19) for $R = r_{\beta n-4}$, $r = r_{\beta n-5}$ and $N = n_{\beta n-4}(a')$, we have that for $\alpha' = (1+2\eta)a'\beta^2/2$ and $n$ large enough,

$$(8.7) \qquad \max_{z \in \mathbb{Z}_{K_n}^2} \mathbf{P}(\mathcal{R}_{\beta n-4}^z(a') > \alpha' t_n^*) \xrightarrow[n \to \infty]{} 0.$$

Further, if $\mathcal{R}_{\beta n-4}^z(a') \leq \alpha' t_n^*$, then $\widetilde{W}^z \geq |\mathcal{L}_{K_n}(\alpha') \cap D(z, r_{\beta n-6})|$. Thus, taking $\eta > 0$ and $\xi > 0$ small enough for $\alpha' < \beta^2$ and $2\beta(1-\alpha) - 4\delta \leq (2\beta - 2\alpha'/\beta)$, we get (8.6) by combining (8.7) and Theorem 1.2.

It thus remains only to show that $p_{n,3} \to 0$ in order to complete the proof of the lemma. To this end, let $\bar{a} = (1 + 2\eta)a$, recall the set $\mathbf{Go}(n)$ [taking $\tilde{n} = n$ in (6.1)] and note that

$$
\begin{aligned}
p_{n,3} \leq{}& K_n^{2\alpha + 2\varepsilon} \max_{(z,x) \in \mathbf{Go}(n)} \mathbf{P}(\mathcal{T}_{K_n}(x) > \mathcal{R}_n^z(\tilde{a}), |\widehat{N}_{n,\beta n}^z(\tilde{a}) - 1| \geq \eta) \\
&+ K_n^{2\alpha + 2\varepsilon} \max_{z \in \mathbb{Z}_{K_n}^2} \mathbf{P}(\widehat{N}_{n,\beta n}^z(\tilde{a}) \leq 1 + \eta, \widehat{N}_{n,\beta n}^z(a) > 1 + 2\eta) \\
&+ K_n^{2\alpha + 2\varepsilon} \max_{z \in \mathbb{Z}_{K_n}^2} \mathbf{P}(\widehat{N}_{\beta n, \beta n-4}^z(\bar{a}) > (1 + 2\eta)^2) \\
:={}& \tilde{p}_n(\beta) + p_{n,4} + p_{n,5},
\end{aligned}
$$

where

$$
\widehat{N}_{\beta n, \beta n-4}^z(\bar{a}) = N_{\beta n, \beta n-4}^z(\bar{a}) / n_{\beta n - 4}(\bar{a})
$$

and the bound above (and in particular the last term $p_{n,5}$) follows from the inclusion

$$
\{N_{n,\beta n-4}^z(\bar{a}) > n_{\beta n-4}(a'), \widehat{N}_{n,\beta n}^z(a) \leq (1 + 2\eta)\} \subset \{\widehat{N}_{\beta n, \beta n-4}^z(\bar{a}) > (1 + 2\eta)^2\},
$$

which is obtained by unraveling the definitions.

Since $\tilde{a} = 2\alpha - \xi$ and $\gamma_1(\beta') = 1$, we have by (6.4) that for any $\beta' \in [\beta(1 - \alpha), \beta]$,

$$
\tag{8.8} \tilde{p}_n(\beta') \leq \sup_{\beta' \in [\beta(1-\alpha), \beta], |\gamma^2 - 1| \geq \eta} K_n^{2\alpha - (2\alpha - \xi)F_{1,\beta'}(\gamma) + 3\varepsilon} \xrightarrow[n \to \infty]{} 0,
$$

for $\varepsilon = \varepsilon(\alpha, \beta, \eta)$ and $\xi = \xi(\alpha, \beta, \eta)$ sufficiently small, using the fact that $(\gamma, \beta') \mapsto F_{1,\beta'}(\gamma)$ is continuous and $F_{1,\beta'}(\gamma) > F_{1,\beta'}(1) = 1$ for $\gamma \neq 1$.

By the strong Markov property of the simple random walk at $\mathcal{R}_n^z(\tilde{a})$ and the bound of (6.3) at $\gamma = ((1 + 2\eta)a - (1 + \eta)\tilde{a})/(a - \tilde{a})$, we have that

$$
\tag{8.9}
\begin{aligned}
p_{n,4} &\leq K_n^{2\alpha + 2\varepsilon} \max_{z \in \mathbb{Z}_{K_n}^2, y \in \partial D(z, r_n)} \mathbf{P}^y(\widehat{N}_{n,\beta n}^z(a - \tilde{a}) \geq \gamma) \\
&\leq K_n^{2\alpha + 3\varepsilon - (a - \tilde{a})F_{0,\beta}(\gamma)} \xrightarrow[n \to \infty]{} 0,
\end{aligned}
$$

for $\xi = \xi(\alpha, \beta, \eta)$ small enough, since $\gamma \to \infty$ and $(a - \tilde{a})F_{0,\beta}(\gamma) = 2\xi F_{0,\beta}(\frac{\eta a}{2\xi} + 1 + \eta) \to \infty$ as $\xi \downarrow 0$.

We complete the proof of the lemma by showing that $p_{n,5} = O(e^{-n^2})$. To this end, first note that by (2.4), the probability that the number of excursions from $\partial D(z, r_{\beta n-5})$ to $\partial D(z, r_{\beta n-4})$ until time $T_{\partial D(z, r_{\beta n})}$ exceeds $2\eta n_{\beta n-4}(\bar{a})$ is bounded for large $n$ and all $z$ by $(9/10)^{\eta n_{\beta n-4}(\bar{a})} = O(e^{-2n^2})$. Hence, using the strong Markov property at $T_{\partial D(z, r_{\beta n})}$ and translation invariance of the simple random walk, it suffices to show that $\mathbf{P}^x(\widehat{N}_{\beta n, \beta n-4}^0(\bar{a}) >$ $1 + 2\eta) = O(e^{-2n^2})$, uniformly in $x \in \partial D(0, r_{\beta n})$. Let $\mathbf{P}_n$ denote probabilities with respect to the random walk in $\mathbb{Z}_{K_n}^2$. Then, uniformly in $x \in \partial D(0, r_{\beta n})$, by conditioning on the $\sigma$-algebra $\mathcal{G}^0$ of excursions from $\partial D(0, r_{\beta n})$ to $\partial D(0, r_{\beta n-1})$ and twice using Lemma 2.4 [for $r = r_{\beta n-5}$, $m = n_{\beta n}(\bar{a})$, first with $K = K_n$ and then with $K = K_{\beta n}$], we see that

$$
\begin{aligned}
(8.10) \qquad & \mathbf{P}_n^x(\widehat{N}_{\beta n, \beta n-4}^0(\bar{a}) > 1 + 2\eta) \\
& \qquad = (1 + o(1_n)) \mathbf{P}_{\beta n}^x(\widehat{N}_{\beta n, \beta n-4}^0(\bar{a}) > 1 + 2\eta).
\end{aligned}
$$

Then, for $\overline{\alpha} = (1 + \eta)\bar{a}/2$, uniformly in $x$ as above $\mathbf{P}_{\beta n}^x(\mathcal{R}_{\beta n}^0(\bar{a}) \geq \overline{\alpha} t_{\beta n}^*) = O(e^{-2n^2})$ by (3.19) and $\mathbf{P}_{\beta n}^x(\mathcal{R}_{\beta n-4}^0((1+2\eta)\bar{a}) \leq \overline{\alpha} t_{\beta n}^*) = O(e^{-2n^2})$ by (3.18). So, the right-hand probability in (8.10) which can be rewritten as $\mathbf{P}_{\beta n}^x(\mathcal{R}_{\beta n}^0(\bar{a}) > \mathcal{R}_{\beta n-4}^0((1+2\eta)\bar{a}))$ is uniformly in $x$ at most $O(e^{-2n^2})$. $\square$

PROOF OF LEMMA 8.2. With $\widehat{\mathcal{Z}}_{n, \beta'}$ as in the proof of Lemma 8.1, let $z_{\beta'}(x)$ denote the point in $\widehat{\mathcal{Z}}_{n, \beta'}$ closest to $x$, and $\mathcal{Z}_{\beta', \eta} = \{z \in \widehat{\mathcal{Z}}_{n, \beta'} : |\widehat{N}_{n, \beta' n}^z(\tilde{a}) - 1| \leq \eta\}$. Taking $h < 2$, to be chosen below, set $\beta_j = \beta(h/2)^j$ for $j = 0, 1, \ldots$ and let $\ell$ be the smallest integer so that $\beta_\ell \leq \beta(1 - \alpha)$. Let $\widehat{W}^x(\cdot, \cdot)$ be as in (8.1), but with the set $\widehat{\mathcal{L}}_{K_n}(\tilde{a})$ of (7.3) instead of $\mathcal{L}_{K_n}(\alpha)$. Note that if $\mathcal{R}_n^z(\tilde{a}) < \alpha t_n^*$ for all $z \in \mathbb{Z}_{K_n}^2$, then $\mathcal{L}_{K_n}(\alpha) \subseteq \widehat{\mathcal{L}}_{K_n}(\tilde{a})$ and $W^x(\cdot, \cdot) \leq \widehat{W}^x(\cdot, \cdot)$. Also, automatically $W^x(0, \beta_\ell) \leq K_n^{2\beta(1-\alpha)}$, so for all $n$ sufficiently large the event $W^x \geq K_n^{2\beta(1-\alpha)+5\delta}$ implies that $W^x(\beta_{j+1}, \beta_j) \geq K_n^{2\beta_j(1-\alpha)+4\delta}$ for some $j = 0, \ldots, \ell - 1$. Thus, we bound the event $\{x \in \mathcal{L}_{K_n}(\alpha), W^x \geq K_n^{2\beta(1-\alpha)+5\delta}\}$ in the definition of $\bar{p}_n$ by the union of the events $\{\mathcal{R}_n^z(\tilde{a}) \geq \alpha t_n^*$ for some $z\}$ and $\{x \in \widehat{\mathcal{L}}_{K_n}(\tilde{a}), \widehat{W}^x(\beta_{j+1}, \beta_j) \geq K_n^{2\beta_j(1-\alpha)+4\delta}\}$, for $j = 0, \ldots, \ell - 1$. Splitting the latter events according to whether $z_{\beta_j}(x) \in \mathcal{Z}_{\beta_j, \eta}$ or not, we get that

$$
\bar{p}_n \leq p_{n,0} + \sum_{j=0}^{\ell-1} \bar{p}_{n,j} + \sum_{j=0}^{\ell-1} \tilde{p}_n(\beta_j),
$$

where

$$
\begin{aligned}
\bar{p}_{n,j} = K_n^{2\alpha+\varepsilon} & \max_{x \in \mathbb{Z}_{K_n}^2} \mathbf{P}(x \in \widehat{\mathcal{L}}_{K_n}(\tilde{a}), z_{\beta_j}(x) \in \mathcal{Z}_{\beta_j, \eta}, \\
& \widehat{W}^x(\beta_{j+1}, \beta_j) \geq K_n^{2\beta_j(1-\alpha)+4\delta}).
\end{aligned}
$$

By (7.4) we know that $p_{n,0} \to 0$ and by (8.8) also $\tilde{p}_n(\beta_j) \to 0$ for $j = 0, \ldots, \ell - 1$.

Turning to deal with $\bar{p}_{n,j}$, let $D_{n,j}(x)$ denote the annulus $D(x, r_{\beta_j n-3}) \setminus D(x, r_{\beta_{j+1} n-3})$. Since, for any $w > 0$ and $x \in \mathbb{Z}_{K_n}^2$,

$$\mathbf{P}(x \in \widehat{\mathcal{L}}_{K_n}(\tilde{a}), z_{\beta_j}(x) \in \mathcal{Z}_{\beta_j, \eta}, \widehat{W}^x(\beta_{j+1}, \beta_j) \geq w)$$
$$\leq w^{-1} \sum_{y \in D_{n,j}(x)} \mathbf{P}(x, y \in \widehat{\mathcal{L}}_{K_n}(\tilde{a}), z_{\beta_j}(x) \in \mathcal{Z}_{\beta_j, \eta}),$$

while $\log r_{\beta_j n-3} / \log K_n \to \beta_j$ and $\gamma_h(\beta_j) \leq \sqrt{1-\eta}$, which we may assume by taking $\eta$ sufficiently small, it follows from (6.5) that for all $n$ large enough

$$\bar{p}_{n,j} \leq K_n^{2\alpha(1+\beta_j)-2\delta} \max_{x,y \in D_{n,j}(x)} \mathbf{P}(x, y \in \widehat{\mathcal{L}}_{K_n}(\tilde{a}), z_{\beta_j}(x) \in \mathcal{Z}_{\beta_j, \eta})$$
$$\leq \sup_{\beta(1-\alpha) \leq \beta' \leq \beta} K_n^{2\alpha(1+\beta')-\tilde{a}F_{h,\beta'}(\sqrt{1-\eta})-\delta}.$$

Then, $\bar{p}_{n,j} \to 0$ as $n \to \infty$ for $\eta, \xi$ sufficiently small and $h < 2$ sufficiently close to 2 using the fact that $(\gamma, h, \beta') \mapsto F_{h,\beta'}(\gamma)$ is continuous and $F_{2,\beta'}(1) = 1 + \beta'$. Possibly decreasing $\varepsilon$ and $\xi$ for (8.8) to hold we complete the proof of (8.3). $\square$

**9. Upper bounds for pairs of late points.** Recall that $F_{h,\beta}(\gamma) = \frac{(1-\gamma\beta)^2}{1-\beta} + h\gamma^2\beta$. We begin by showing that

$$(9.1) \qquad 2 + 2\beta - 2\alpha \inf_{\gamma \in \Gamma_{\alpha,\beta}} F_{2,\beta}(\gamma)$$
$$= \begin{cases} 2 + 2\beta - 4\alpha/(2-\beta), & \text{if } \beta \leq 2(1-\sqrt{\alpha}), \\ 8(1-\sqrt{\alpha}) - 4(1-\sqrt{\alpha})^2/\beta, & \text{if } \beta \geq 2(1-\sqrt{\alpha}), \end{cases}$$

where $\Gamma_{\alpha,\beta} = \{\gamma \geq 0 : 2 - 2\beta - 2\alpha F_{0,\beta}(\gamma) \geq 0\}$, thereby establishing the equivalence of (1.8) and (1.11). Indeed, as noted before, $F_{2,\beta}(\gamma)$ is quadratic, with minimum value $F_{2,\beta}(\gamma_2) = 2/(2-\beta)$ achieved at $\gamma_2(\beta) = 1/(2-\beta) < 1$. It is easy to check that $\Gamma_{\alpha,\beta}$ is the interval $[\gamma_-, \gamma_+]$ for

$$(9.2) \qquad \gamma_{\pm} = \gamma_{\pm}(\alpha, \beta) = \beta^{-1} \max\{1 \pm \alpha^{-1/2}(1-\beta), 0\}.$$

Since $\gamma_2 < 1 < \gamma_+$ we see that $\gamma_2 \in \Gamma_{\alpha,\beta}$ if and only if $\gamma_- \leq \gamma_2$, leading to the explicit formula

$$(9.3) \qquad \rho(\alpha, \beta) = 2 + 2\beta - 2\alpha F_{2,\beta}(\max\{\gamma_-, \gamma_2\})$$

[where we denote hereafter the left-hand side of (9.1) as $\rho(\alpha, \beta)$]. Combining this with the fact that $\gamma_-(\alpha, \beta) > \gamma_2(\beta)$ is equivalent to $\beta > 2(1-\sqrt{\alpha})$, we obtain the identity (9.1). Clearly, $\beta \mapsto \rho(\alpha, \beta)$ is continuous on $(0, 1)$ and by (9.1) it is also monotone increasing in $\beta$ [for $\beta \geq 2(1-\sqrt{\alpha})$ by inspection, while for $\beta \leq 2(1-\sqrt{\alpha})$ we have that $d\rho/d\beta \geq 1$].

We prove in this section that for any $0 < \alpha, \beta, \delta < 1$,

$$(9.4) \qquad \lim_{K \to \infty} \mathbf{P}(|\{(x,y) : x, y \in \mathcal{L}_K(\alpha), d(x,y) \leq K^\beta\}| \geq K^{\rho(\alpha,\beta) + 4\delta}) = 0.$$

To this end, let

$$(9.5) \qquad \Psi_{\alpha, \beta_2, \beta_1, n} = \{(x,y) : x, y \in \mathcal{L}_{K_n}(\alpha), r_{\beta_2 n-3} < d(x,y) \leq r_{\beta_1 n-3}\}.$$

It suffices (as usual) to prove that (9.4) holds for $K_n = n^{\bar\gamma}(n!)^3$, uniformly in $\bar\gamma \in \mathcal{I}$. Further, $\log r_{\beta n-3} / \log K_n \to \beta$, so fixing $0 < \alpha, \beta, \delta < 1$, it is enough to show that

$$(9.6) \qquad \lim_{n \to \infty} \mathbf{P}(|\Psi_{\alpha, 0, \beta, n}| \geq K_n^{\rho(\alpha,\beta) + 4\delta}) = 0.$$

Note that $|\Psi_{\alpha, 0, \beta(1-\alpha), n}| \leq K_n^{2\beta(1-\alpha)} |\mathcal{L}_{K_n}(\alpha)|$ for some universal $n_0 = n_0(\alpha, \beta) < \infty$ and all $n \geq n_0$, while

$$\rho(\alpha, \beta) \geq 2 + 2\beta - 2\alpha F_{2,\beta}(1) = 2(1-\alpha) + 2\beta(1-\alpha),$$

so that it follows from (1.2) that

$$\begin{aligned}(9.7) \qquad &\lim_{n \to \infty} \mathbf{P}(|\Psi_{\alpha, 0, \beta(1-\alpha), n}| \geq K_n^{\rho(\alpha,\beta) + 4\delta}) \\ &\qquad \leq \lim_{n \to \infty} \mathbf{P}(|\mathcal{L}_{K_n}(\alpha)| \geq K_n^{2(1-\alpha) + 4\delta}) = 0.\end{aligned}$$

The following lemma will be proven below.

LEMMA 9.1. *We can choose $h < 2$ sufficiently close to 2 and $\tilde{a} < 2\alpha$ sufficiently close to $2\alpha$ such that for any $\beta' \in [\beta(1-\alpha), \beta]$*

$$(9.8) \qquad q_{n,\beta'} := \mathbf{P}(|\widehat\Psi_{\tilde{a}, h, \beta', n}| \geq K_n^{\rho(\alpha,\beta') + 3\delta}) \underset{n \to \infty}{\longrightarrow} 0,$$

*where*

$$\widehat\Psi_{\tilde{a}, h, \beta', n} = \{(x,y) : x, y \in \widehat{\mathcal{L}}_{K_n}(\tilde{a}), r_{\beta' h n/2-3} < d(x,y) \leq r_{\beta' n-3}\}.$$

Fix $h < 2$, $\tilde{a} < 2\alpha$ according to Lemma 9.1. We then set $\beta_j = \beta(h/2)^j$ and $\ell$ as the smallest integer such that $\beta_\ell \leq \beta(1-\alpha)$. By Lemma 9.1 we have that $q_{n, \beta_j} \longrightarrow 0$, as $n \to \infty$ for $j = 0, \ldots, \ell-1$. Combining this with (7.5), the monotonicity of $\beta \mapsto \rho(\alpha, \beta)$ and (9.7), we establish (9.6).

PROOF OF LEMMA 9.1. Let $D_{n,\beta'}(x)$ denote the annulus $D(x, r_{\beta' n-3}) \setminus D(x, r_{\beta' h n/2-3})$. Fix $0 < \eta < 1$ to be chosen below, abbreviating $\gamma_- = \gamma_-(\alpha, \beta')$, $\gamma_* = (1-\eta)\gamma_-(\tilde{a}/2, \beta')$ and $\gamma_h = \gamma_h(\beta')$. We will argue separately depending on whether or not $\gamma_* \leq \gamma_h$. Consider first the case where $\gamma_* \leq \gamma_h$. Applying (6.5) at $\gamma = \gamma_h$ we conclude that for all $n$ large enough,

$$\max_{x \in \mathbb{Z}_{K_n}^2} \max_{y \in D_{n,\beta'}(x)} \mathbf{P}(x, y \in \widehat{\mathcal{L}}_{K_n}(\tilde{a})) \leq K_n^{-\tilde{a} F_{h,\beta'}(\gamma_h) + \delta}.$$

By (9.3) at $\beta'$, this implies that if $\gamma_* \leq \gamma_h$, then

$$q_{n,\beta'} \leq K_n^{-\rho(\alpha,\beta')-3\delta} \sum_{x \in \mathbb{Z}_{K_n}^2} \sum_{y \in D_{n,\beta'}(x)} \mathbf{P}(x,y \in \widehat{\mathcal{L}}_{K_n}(\tilde{a})) \leq K_n^{\bar{g}_{\beta'}(\eta,\tilde{a},h)-\delta},$$

where $\bar{g}_{\beta'}(\eta,\tilde{a},h) = 2\alpha F_{2,\beta'}(\max\{\gamma_-,\gamma_2\}) - \tilde{a}F_{h,\beta'}(\max\{\gamma_*,\gamma_h\})$. (Here $\max\{\gamma_*,\gamma_h\} = \gamma_h$.) Note that $\bar{g}_{\beta'}(0,2\alpha,2) = 0$ for all $\beta'$; hence for any $\delta > 0$ we can and shall take $h$ sufficiently close to 2, $\tilde{a} < 2\alpha$ sufficiently close to $2\alpha$ and $\eta > 0$ sufficiently small so that $\bar{g}_{\beta'}(\eta,\tilde{a},h) < \delta/2$ for all $\beta' \in [\beta(1-\alpha),\beta]$. Clearly, this choice of parameters guarantees that $q_{n,\beta'} \underset{n\to\infty}{\longrightarrow} 0$ whenever $\gamma_* \leq \gamma_h$.

Keeping this choice of $h$, $\tilde{a}$ and $\eta$, we turn to deal with the case where $\gamma_* > \gamma_h$, denoting by $\widehat{\mathcal{Z}}_{n,\beta'}$ the subgrid in $\mathbb{Z}_{K_n}^2$ of spacing $4r_{\beta'n-4}$. Let $z_{\beta'}(x)$ denote the point closest to $x$ in $\widehat{\mathcal{Z}}_{n,\beta'}$ so

$$\mathbf{P}(\widehat{N}_{n,\beta'n}^{z_{\beta'}(x)}(\tilde{a}) \leq \gamma_*^2) \leq \mathbf{P}\left( \min_{z \in \widehat{\mathcal{Z}}_{n,\beta'}} \{\widehat{N}_{n,\beta'n}^z(\tilde{a})\} \leq \gamma_*^2 \right) =: q_n(\beta').$$

Then, using again (9.3) at $\beta'$ and the bound (6.5), now for $\gamma = \gamma_* \geq \gamma_h$, we get that

$$q_{n,\beta'} \leq q_n(\beta') + K_n^{-\rho(\alpha,\beta')-3\delta} \sum_{x \in \mathbb{Z}_{K_n}^2} \sum_{y \in D_{n,\beta'}(x)} \mathbf{P}(x,y \in \widehat{\mathcal{L}}_{K_n}(\tilde{a}), \widehat{N}_{n,\beta'n}^{z_{\beta'}(x)}(\tilde{a}) \geq \gamma_*^2)$$

$$\leq q_n(\beta') + K_n^{2\alpha F_{2,\beta'}(\max\{\gamma_-,\gamma_2\})-2\delta}$$
$$\times \max_{z \in \mathbb{Z}_{K_n}^2} \max_{\substack{x,y \in D(z,r_{\beta'n-2}) \\ d(x,y) \geq r_{\beta'hn/2-3}}} \mathbf{P}(x,y \in \widehat{\mathcal{L}}_{K_n}(\tilde{a}), \widehat{N}_{n,\beta'n}^z(\tilde{a}) \geq \gamma_*^2)$$

$$\leq q_n(\beta') + K_n^{\bar{g}_{\beta'}(\eta,\tilde{a},h)-\delta}.$$

(Here $\max\{\gamma_*,\gamma_h\} = \gamma_*$.) As we have seen, our choice of parameters guarantees that $\bar{g}_{\beta'}(\eta,\tilde{a},h) < \delta/2$. Moreover, since $\gamma_* \leq 1 \leq \gamma_0$ [for any $\beta' \in (0,1)$], it follows by (6.3) that for any $\varepsilon > 0$ and all $n$ large enough,

$$(9.9) \qquad q_n(\beta') \leq |\widehat{\mathcal{Z}}_{n,\beta'}| K_n^{-\tilde{a}F_{0,\beta'}(\gamma_*)+\varepsilon} \leq K_n^{2-2\beta'-\tilde{a}F_{0,\beta'}(\gamma_*)+2\varepsilon}.$$

Note that for $h \leq 2$ we have $\gamma_h \geq 1/2$. Hence, using our assumption that $\gamma_* \geq \gamma_h$ and the definition of $\gamma_*$, we have that $\gamma_-(\tilde{a}/2,\beta') \geq 1/2 > 0$. This guarantees that $\gamma_-(\tilde{a}/2,\beta')$ is the lower boundary of $\{\gamma: 2-2\beta'-\tilde{a}F_{0,\beta'}(\gamma) \geq 0\}$. It follows that $2-2\beta'-\tilde{a}F_{0,\beta'}(\gamma_*) < 0$ uniformly in $\beta' \in [\beta(1-\alpha),\beta]$ for which $\gamma_* \geq \gamma_h$. Hence we can find $\varepsilon > 0$ so that $q_n(\beta') \underset{n\to\infty}{\longrightarrow} 0$ uniformly in this set of values of $\beta'$, implying in turn that $q_{n,\beta'} \underset{n\to\infty}{\longrightarrow} 0$. This completes the proof of Lemma 9.1. $\square$

**10. Lower bounds for pairs of late points.** Fix $0 < \alpha, \beta < 1$. Recall the notation $K_n = n^{\bar{\gamma}}(n!)^3$ and the sets $\Psi_{\alpha,0,\beta,n}$ of (9.5). We show that if $\gamma_-(\alpha,\beta) < \gamma < 1$ and $1 - \alpha > \delta > \xi > 0$ are such that $2 - 2\beta - (2\alpha + \xi)F_{2,\beta}(\gamma) > 2\delta$, then

$$(10.1) \qquad \lim_{n \to \infty} \mathbf{P}(|\Psi_{\alpha,0,\beta,n}| \geq K_n^{2+2\beta-(2\alpha+\xi)F_{2,\beta}(\gamma)-5\delta}) = 1,$$

uniformly in $\bar{\gamma} \in \mathcal{I}$. In view of (9.3), taking $\xi, \delta \downarrow 0$ followed by $\gamma \in (\gamma_-(\alpha,\beta), 1)$ that converges to $\max(\gamma_-(\alpha,\beta), \gamma_2(\beta))$, we get the lower bound in Theorem 1.4 for the subsequence $K_n$. By the uniformity in $\bar{\gamma}$ this bound extends to all integers.

Fixing $\gamma$, $\delta$ and $\xi$ as above, set $a = 2\alpha + \xi$, recall the notation $r_k$, $n_k(a)$, $\mathcal{R}_k^x(a)$ and $N_{k,l}^x(a)$ of Section 4 and let

$$(10.2) \qquad \hat{n}_k = 3a^*\left(k - \frac{(\beta - \gamma\beta)}{(1 - \gamma\beta)}n\right)^2 \log k, \qquad \beta n \leq k \leq n,$$

where $a^* = a(1 - \gamma\beta)^2/(1-\beta)^2$, so that $\hat{n}_n = n_n(a)$ and $\hat{n}_{\beta n} = \gamma^2 n_{\beta n}(a)$.

Let $\mathcal{Z}_n \subset \mathbb{Z}_{K_n}^2$ be a maximal set of points in $\mathbb{Z}_{K_n}^2 \setminus D(0, r_n)$ which are $4r_{\beta n+4}$ separated, such that $(0, 2r_n) \in \mathcal{Z}_n$. We will say that a point $z \in \mathcal{Z}_n$ is $(n, \beta)$-qualified if $N_{n,k}^z \overset{k}{\sim} \hat{n}_k$ for all $\beta n \leq k \leq n-1$ and in addition

$$\widetilde{W}^z := |\{y \in D(z, r_{\beta n-4}) : \mathcal{T}_{K_n}(y) > \mathcal{R}_n^z(a)\}| \geq K_n^{\beta(2-a\gamma^2)-2\delta}$$

(compare with the definition of $n$-successful points in Section 4). If

$$\min_{z \in \mathbb{Z}_{K_n}^2} \mathcal{R}_n^z(a) \geq \alpha t_n^*,$$

then

$$|\Psi_{\alpha,0,\beta,n}| \geq \sum_{z \in \mathcal{Z}_n} (\widetilde{W}^z)^2 \geq |\{z \in \mathcal{Z}_n : z \text{ is } (n,\beta)\text{-qualified}\}| K_n^{2\beta(2-a\gamma^2)-4\delta}.$$

Since $\mathbf{P}(\min_{z \in \mathbb{Z}_{K_n}^2} \mathcal{R}_n^z(a) \leq \alpha t_n^*) \to 0$ as $n \to \infty$ (see the term $p_{n,1}$ in the proof of Lemma 8.1), and $(1 - \beta)a^* = aF_{0,\beta}(\gamma)$, we thus get (10.1) as soon as we show that

$$(10.3) \qquad \lim_{n \to \infty} \mathbf{P}(|\{z \in \mathcal{Z}_n : z \text{ is } (n,\beta)\text{-qualified}\}| \geq K_n^{(1-\beta)(2-a^*)-\delta}) = 1.$$

The following analogue of Lemma 4.2, whose proof is deferred to the end of this section, is the key to the proof of (10.3).

LEMMA 10.1. *For any $x, y \in \mathcal{Z}_n$, let $l(x,y) = \max\{k : D(x, r_k+1) \cap D(y, r_k+1) = \varnothing\} \wedge n$ [note that $l(x,y) \geq \beta n + 4$]. Then there exist $b \geq 10$ and $\hat{q}_n \geq (r_n/r_{\beta n})^{-a^*+o(1_n)}$ such that*

$$(10.4) \qquad \mathbf{P}(z \text{ is } (n,\beta)\text{-qualified}) = (1 + o(1_n))\hat{q}_n,$$

*uniformly in $\bar{\gamma} \in \mathcal{I}$ and $z \in \mathcal{Z}_n$. Furthermore, for any $\varepsilon > 0$ we can find $C = C(b, \varepsilon) < \infty$ such that for all $n$ and any $x, y \in \mathcal{Z}_n$ with $l(x, y) < n$,*

$$(10.5) \quad \mathbf{P}(x, y \text{ are both } (n, \beta)\text{-qualified}) \leq \widehat{q}_n^{~2} C^{n-l(x,y)} n^b \left( \frac{r_n}{r_{l(x,y)}} \right)^{a^*+\varepsilon},$$

*while for all $n$ and $x, y \in \mathcal{Z}_n$ with $l(x, y) = n$,*

$$(10.6) \qquad \mathbf{P}(x, y \text{ are both } (n, \beta)\text{-qualified}) \leq (1 + o(1_n))\widehat{q}_n^{~2}.$$

The proof of (10.3) then proceeds exactly as the proof of (4.3), where the condition $2 - 2\beta - aF_{0,\beta}(\gamma) > 2\delta$ implies that $a^* < 2$ and by (10.4) the expected number of $(n, \beta)$-qualified points is $K_n^{(1-\beta)(2-a^*)+o(1_n)}$. So, with

$$V_\ell = \sum_{x,y \in \mathcal{Z}_n, l(x,y)=\ell} \mathbf{P}(x, y \text{ are both } (n, \beta)\text{-qualified}), \qquad \ell = \beta n + 4, \ldots, n,$$

it suffices by (10.6) to show that

$$(10.7) \qquad \qquad \sum_{\ell=\beta n+4}^{n-1} V_\ell \leq o(1_n) |\mathcal{Z}_n|^2 \widehat{q}_n^{~2}.$$

With $C_m$ denoting generic finite constants that are independent of $n$, for any $\ell \in [\beta n + 4, n)$ and $x \in \mathcal{Z}_n$ there are at most $C_0 r_{\ell+1}^2 / r_{\beta n+4}^2$ points $y \in \mathcal{Z}_n \cap D(x, 2(r_{\ell+1} + 1))$. Consequently, we have by (10.5) and the definition of $l(x, y)$ that for any $\ell \in [\beta n + 4, n)$,

$$V_\ell \leq C_2 |\mathcal{Z}_n| \left( \frac{r_{\ell+1}}{r_{\beta n+4}} \right)^2 \widehat{q}_n^{~2} n^b C^{n-\ell} \left( \frac{r_n}{r_\ell} \right)^{a^*+\varepsilon}.$$

Similarly to the derivation of (4.9), taking $\varepsilon < 2 - a^*$ and summing over $\ell$ results with (10.7), hence completing the proof of (10.3).

PROOF OF LEMMA 10.1. Let $\mathcal{R}_{\beta n,m}^z$ denote the time until completion of the first $m$ excursions from $\partial D(z, r_{\beta n-1})$ to $\partial D(z, r_{\beta n})$, and set $\widehat{\mathbf{A}}_m^z = \{\widetilde{W}_m^z \geq K_n^{\beta(2-a\gamma^2)-2\delta}\}$ for $\widetilde{W}_m^z = |\{y \in D(z, r_{\beta n-4}) : \mathcal{T}_{K_n}(y) > \mathcal{R}_{\beta n,m}^z\}|$. Recall that $\widehat{n}_{\beta n} = n_{\beta n}(\gamma^2(2\alpha + \xi))$, so applying (3.19) with $R = r_{\beta n}$, $r = r_{\beta n-1}$ and $N = \widehat{n}_{\beta n} + \beta n$, we see that for all $m \leq \widehat{n}_{\beta n} + \beta n$,

$$\mathbf{P}(\widetilde{W}_m^z \geq |D(z, r_{\beta n-4}) \cap \mathcal{L}_{K_n}((\alpha + \xi)\gamma^2 \beta^2)|)$$
$$\geq \mathbf{P}(\mathcal{R}_{\beta n,m}^z \leq (\alpha + \xi)\gamma^2 \beta^2 t_n^*) = 1 - o(1_n).$$

Hence, by Theorem 1.2 we have that

$$(10.8) \qquad \mathbf{P}(\widehat{\mathbf{A}}_m^z) = 1 - o(1_n) \qquad \text{uniformly in } m \overset{\beta n}{\sim} \widehat{n}_{\beta n}.$$

Starting at $0 \notin D(z, r_{\beta n})$ we see that the event $\widehat{\mathbf{A}}_m^z$ belongs to the $\sigma$-algebra $\mathcal{H}^z(m)$ corresponding to $r = r_{\beta n-2}$, $R = r_{\beta n-1}$ and $R' = r_{\beta n}$ in Lemma 2.4. Further, if the event $\{N_{n,\beta n}^z = m\} \in \mathcal{G}_{\beta n}^z$ occurs, then the law of $\mathbf{A}_{\beta n}^z = \{\widetilde{W}^z \geq K_n^{\beta(2-a\gamma^2)-2\delta}\}$ conditioned upon $\mathcal{G}_{\beta n}^z$ is the same as the law of $\widehat{\mathbf{A}}_m^z$ conditioned upon $\mathcal{G}_{\beta n}^z$. Consequently, by Lemma 2.4 and (10.8), uniformly in $m \stackrel{\beta n}{\sim} \widehat{n}_{\beta n}$,

$$(10.9) \quad \mathbf{P}(\mathbf{A}_{\beta n}^z | \mathcal{G}_{\beta n}^z, N_{n,\beta n}^z = m) = \mathbf{P}(\widehat{\mathbf{A}}_m^z | \mathcal{G}_{\beta n}^z, N_{n,\beta n}^z = m) = 1 - o(1_n).$$

With $M_l = \{l, \ldots, n-1\}$, by (10.9) and the fact that $\{N_{n,k}^z \stackrel{k}{\sim} \widehat{n}_k, k \in M_{\beta n}\} \in \mathcal{G}_{\beta n}^z$ we get that

$$\mathbf{P}(z \text{ is } (n, \beta)\text{-qualified}) = \mathbf{P}(N_{n,k}^z \stackrel{k}{\sim} \widehat{n}_k, k \in M_{\beta n}; \mathbf{A}_{\beta n}^z)$$

$$= \sum_{m \stackrel{\beta n}{\sim} \widehat{n}_{\beta n}} \mathbb{E}(N_{n,k}^z \stackrel{k}{\sim} \widehat{n}_k, k \in M_{\beta n+1};$$

$$N_{n,\beta n}^z = m; \mathbf{P}(\mathbf{A}_{\beta n}^z | \mathcal{G}_{\beta n}^z, N_{n,\beta n}^z = m))$$

$$= (1 + o(1_n))\mathbf{P}(N_{n,k}^z \stackrel{k}{\sim} \widehat{n}_k, k \in M_{\beta n}).$$

Therefore, taking $m_n = \widehat{n}_n = n_n(a)$, by (5.9) we get (10.4) for

$$(10.10) \quad \widehat{q}_n = \sum_{\substack{m_{\beta n}, \ldots, m_{n-1} \\ |m_\ell - \widehat{n}_\ell| \leq \ell}} \prod_{\ell=\beta n}^{n-1} \binom{m_{\ell+1} + m_\ell - 1}{m_\ell} p_\ell^{m_\ell} (1 - p_\ell)^{m_{\ell+1}}.$$

It is not hard to check that our choice (10.2) implies that for some $C < \infty$ and all $k \in M_{\beta n}$, if $|m - \widehat{n}_k| \leq k$ and $|l + 1 - \widehat{n}_{k+1}| \leq k+1$, then

$$\left| \frac{m}{l} - 1 - \frac{2}{k - ((\beta - \gamma\beta)/(1 - \gamma\beta))n} \right| \leq \frac{C}{k \log k},$$

which by adapting the proof of [3], Lemma 7.2, shows that uniformly in $m_k \stackrel{k}{\sim} \widehat{n}_k$ and $m_{k+1} \stackrel{k+1}{\sim} \widehat{n}_{k+1}$,

$$(10.11) \quad \frac{C'k^{-3a^*-1}}{\sqrt{\log k}} \leq \binom{m_{k+1} + m_k - 1}{m_k} p_k^{m_k} (1 - p_k)^{m_{k+1}} \leq \frac{Ck^{-3a^*-1}}{\sqrt{\log k}}$$

with $0 < C', C < \infty$ independent of $k$. Putting (10.10) and (10.11) together we see that $\widehat{q}_n = (r_n/r_{\beta n})^{-a^* + o(1_n)}$ as claimed.

It suffices to prove the upper bounds of (10.5) and (10.6) with the events $\{z \text{ is } (n, \beta)\text{-qualified}\}$ replaced by the larger events $\mathcal{A}(z, n, \beta) := \{N_{n,k}^z \stackrel{k}{\sim} \widehat{n}_k, k \in M_{\beta n}\}$. The proof is a rerun of the argument used in Section 5 to prove (4.5) and (4.6), respectively, replacing the events $\{z \text{ is } n\text{-successful}\}$ by $\mathcal{A}(z, n, \beta)$, taking $\rho = \beta$ and $\beta n + 4$ instead of $\rho' n$, excluding 0 from the sets $J_l$ and $I_l$ and replacing everywhere there $\bar{q}_n$ with $\hat{q}_n$, $n_k$ with $\hat{n}_k$ and $a$ with $a^*$. Indeed, the effect of the values $\hat{n}_k$ is in the application of (10.11) whenever (5.7) is used in Section 5. $\square$

**11. Complements and unsolved problems.** (A) Let $L_n^x$ denote the number of times that $x \in \mathbb{Z}^2$ is visited by the simple random walk in $\mathbb{Z}^2$ up to the time $T_{\partial D(0,n)}$ of exit from the disc of radius $n$. For any $0 < \alpha < 1$, set

$$(11.1) \qquad \Psi_n(\alpha) = \left\{ x \in D(0, n) : \frac{L_n^x}{(\log n)^2} \geq 4\alpha/\pi \right\}.$$

Since $\log T_{\partial D(0,n)}/\log n \to 2$ almost surely as $n \to \infty$ (see, e.g., [8], equation (6)), our result ([3], (1.3)) is equivalent to

$$(11.2) \qquad \lim_{n \to \infty} \frac{\log |\Psi_n(\alpha)|}{\log n} = 2(1 - \alpha) \qquad \text{a.s.}$$

Following the line of reasoning of this paper, we expect that for any $0 < \alpha, \beta < 1$, choosing $Y_n$ uniformly in $\Psi_n(\alpha)$,

$$(11.3) \qquad \limsup_{n \to \infty} \frac{\log |\Psi_n(\alpha) \cap D(Y_n, n^\beta)|}{\log n} = 2\beta(1 - \alpha) \qquad \text{a.s.}$$

We also expect that the analysis in this paper can be extended to yield

$$(11.4) \qquad \lim_{n \to \infty} \frac{\log |\{x, y \in \Psi_n(\alpha) : d(x, y) \leq n^\beta\}|}{\log n} = \rho(\alpha, \beta) \qquad \text{a.s.}$$

(B) Our study of planar random walk suggests that the analogous results hold for the planar Wiener sausage. Let $S_\varepsilon(t) = \{x \in \mathbb{T}^2 : \exists\, s \leq t, |W_s - x| \leq \varepsilon\}$ denote the set covered by the Wiener sausage up to time $t$, where $W_t$ is the Brownian motion on the two-dimensional torus $\mathbb{T}^2$. Consider the uncovered set $U_\varepsilon(\alpha) = \mathbb{T}^2 \setminus S_\varepsilon(2\alpha(\log \varepsilon)^2/\pi)$ for $0 < \alpha < 1$ (in [4] we show that $U_\varepsilon(\alpha)$ is empty if $\alpha > 1$). With $\mathcal{L}eb$ denoting Lebesgue's measure, we then expect that

$$(11.5) \qquad \lim_{\varepsilon \to 0} \frac{\log \mathcal{L}eb(U_\varepsilon(\alpha))}{\log \varepsilon} = 2\alpha \qquad \text{a.s.}$$

and for any $x \in \mathbb{T}^2$, $1 > \beta > \sqrt{\alpha}$,

$$(11.6) \qquad \lim_{\varepsilon \to 0} \frac{\log \mathcal{L}eb(U_\varepsilon(\alpha) \cap D(x, \varepsilon^{1-\beta}))}{\log \varepsilon} = 2 - (2\beta - 2\alpha/\beta) \qquad \text{a.s.}$$

We also expect that for $0 < \alpha, \beta < 1$ and $Y_\varepsilon$ chosen according to Lebesgue measure on $U_\varepsilon(\alpha)$,

$$(11.7) \qquad \lim_{\varepsilon \to 0} \frac{\log \mathcal{L}eb(U_\varepsilon(\alpha) \cap D(Y_\varepsilon, \varepsilon^{1-\beta}))}{\log \varepsilon} = 2 - 2\beta(1 - \alpha) \qquad \text{a.s.}$$

and that

$$(11.8) \quad \lim_{\varepsilon \to 0} \frac{\log \int_{U_\varepsilon(\alpha)} \mathcal{L}eb(U_\varepsilon(\alpha) \cap D(x, \varepsilon^{1-\beta}))\, dx}{\log \varepsilon} = 4 - \rho(\alpha, \beta) \qquad \text{a.s.}$$

We believe that these results can be derived by arguments similar to those used here, but have not verified it.

**Acknowledgment.** We thank Bertrand Duplantier for pointing our attention to the work of Brummelhuis and Hilhorst [1] in connection with our earlier work [4] on cover times.

A. DEMBO
DEPARTMENTS OF MATHEMATICS
AND STATISTICS
STANFORD UNIVERSITY
STANFORD, CALIFORNIA 94305
USA
E-MAIL: amir@math.stanford.edu

J. ROSEN
DEPARTMENT OF MATHEMATICS
CITY UNIVERSITY OF NEW YORK
COLLEGE OF STATEN ISLAND
STATEN ISLAND, NEW YORK 10314
USA
E-MAIL: jrosen3@earthlink.net

Y. PERES
DEPARTMENTS OF MATHEMATICS
AND STATISTICS
UNIVERSITY OF CALIFORNIA
BERKELEY, CALIFORNIA 94720
USA
E-MAIL: peres@stat.berkeley.edu

O. ZEITOUNI
DEPARTMENTS OF ELECTRICAL ENGINEERING
AND MATHEMATICS
TECHNION, HAIFA 32000
ISRAEL
AND
DEPARTMENT OF MATHEMATICS
UNIVERSITY OF MINNESOTA
MINNEAPOLIS, MINNESOTA 55455
USA
E-MAIL: zeitouni@math.umn.edu